\documentclass [11pt]{amsart}
\usepackage {etex, amsmath, amssymb, amscd, mathrsfs, hyperref, enumerate, color, slashed, url}
\usepackage[text={6in,8.5in},centering,letterpaper,dvips]{geometry}
\usepackage[all,cmtip]{xy}

\usepackage{mathabx} 
\usepackage{amsthm}
\usepackage{verbatim} 
\usepackage[pdftex]{graphicx}
\usepackage{epstopdf}
\usepackage[export]{adjustbox}
\usepackage{caption, subcaption} 
\usepackage{makecell}
\usepackage{tikz-cd}
\usetikzlibrary{arrows}
\DeclareGraphicsExtensions{.eps,.pdf}

\newtheorem {theorem}{Theorem}[section]
\newtheorem {lemma}[theorem]{Lemma}

\newtheorem {corollary}[theorem]{Corollary}

\newtheorem {definition}[theorem]{Definition}

\newtheorem{exercise}[theorem]{Exercise}

\theoremstyle{remark}
\newtheorem {remark}[theorem]{Remark}
\newtheorem {example}[theorem]{Example}

\def \partiali {\partial\!\!\!/\,}
\def \Di {D\!\!\!\!/\,}
\DeclareFontFamily{U}{mathx}{\hyphenchar\font45}
\DeclareFontShape{U}{mathx}{m}{n}{
      <5> <6> <7> <8> <9> <10>
      <10.95> <12> <14.4> <17.28> <20.74> <24.88>
      mathx10
      }{}
\DeclareSymbolFont{mathx}{U}{mathx}{m}{n}
\DeclareFontSubstitution{U}{mathx}{m}{n}
\DeclareMathAccent{\widecheck}{0}{mathx}{"71}

\def\polhk#1{\setbox0=\hbox{#1}{\ooalign{\hidewidth
    \lower1.5ex\hbox{`}\hidewidth\crcr\unhbox0}}}  
\def\Z {{\mathbb{Z}}}
\def\R {{\mathbb{R}}}

\def\Q {{\mathbb{Q}}}

\def\Conf {\mathcal{C}}

\def\F {\mathbb{F}}
\def\del {\partial}

\def\swf{\mathit{SWF}}
\def\swfh{\mathit{SWFH}}

\def\pin {\operatorname{Pin}(2)}

\newcommand{\s}{\mathfrak{s}}
\newcommand{\bunderline}[1]{\underline{#1\mkern-2mu}\mkern2mu }
\newcommand{\sunderline}[1]{\underline{#1\mkern-3mu}\mkern3mu }
\def\du {\bar{d}}
\def\dl {\bunderline{d}}

\def\Vl {\sunderline{V}}

\begin{document}

\title[Lectures on the triangulation conjecture]{Lectures on the triangulation conjecture}
\author[Ciprian Manolescu]{Ciprian Manolescu}
\address {Department of Mathematics, UCLA, 520 Portola Plaza\\ Los Angeles, CA 90095}
\email {cm@math.ucla.edu}

\begin{abstract}
We outline the proof that non-triangulable manifolds exist in any dimension greater than four. The arguments involve homology cobordism invariants coming from the $\pin$ symmetry of the Seiberg-Witten equations. We also explore a related construction, of an involutive version of Heegaard Floer homology.
\end {abstract}

\maketitle

\section{Introduction}

The triangulation conjecture stated that every topological manifold can be triangulated. The work of Casson \cite{am90} in the 1980's provided counterexamples in dimension $4$. The main purpose of these notes is to describe the proof of the following theorem.

\begin{theorem}[\cite{man13}]
\label{thm:tc}
There exist non-triangulable $n$-dimensional topological manifolds for every $n \geq 5$.
\end{theorem}

The proof relies on previous work by Galewski-Stern \cite{gs80} and Matumoto \cite{mat78}, who reduced this problem to a different one, about the homology cobordism group in three dimensions. Homology cobordism can be explored using the techniques of gauge theory, as was done, for example, by Fintushel and Stern \cite{fs85, fs90}, Furuta \cite{fur90}, and Fr{\o}yshov \cite{fro10}. In 
 \cite{man13}, $\pin$-equivariant Seiberg-Witten Floer homology is used to construct three new invariants of homology cobordism, called $\alpha$, $\beta$ and $\gamma$. The properties of $\beta$ suffice to answer the question raised by Galewski-Stern and Matumoto, and thus prove Theorem~\ref{thm:tc}.

The paper is organized as follows.

Section~\ref{sec:triangulations} contains background material about triangulating manifolds. In particular, we sketch the arguments of Galewski-Stern and Matumoto that reduced Theorem~\ref{thm:tc} to a problem about homology cobordism.

In Section~\ref{sec:sw} we introduce the Seiberg-Witten equations, finite dimensional approximation, and the Conley index. Using these ingredients, we review the construction of  Seiberg-Witten Floer stable homotopy types, following \cite{man03}.

In Section~\ref{sec:abc} we explore the module structure on Borel homology, and more specifically on the $\pin$-equivariant homology of the Seiberg-Witten Floer stable homotopy type. Using this module structure, we define the three numerical invariants $\alpha, \beta,\gamma$, and show that they are preserved by homology cobordism.

Section~\ref{sec:duality} contains material about equivariant Spanier-Whitehead duality. This is applied to understanding the behavior of $\alpha, \beta,\gamma$ under orientation reversal. Showing that $\beta(-Y) = -\beta(Y)$ completes the proof of Theorem~\ref{thm:tc}.

In Section~\ref{sec:hfi} we outline the construction of involutive Heegaard Floer homology, joint work of Hendricks and the author \cite{hm15}. Involutive Heegaard Floer homology is a more computable counterpart to $\Z/4$-equivariant Seiberg-Witten Floer homology, and has its own applications to questions about homology cobordism. 

\medskip
 \textbf {Acknowledgements.} 
These are notes from a set of five lectures delivered by the author at the 22nd G{\"o}kova Geometry / Topology Conference, May 25-30, 2015. I would like to heartfully thank Eylem Zeliha Y{\i}ld{\i}z for taking and preparing the notes, Selman Akbulut for his interest in this work, and Larry Taylor for  correcting an error in a previous version. Partial support came from the NSF grant DMS-1402914.

\section{Triangulations}
\label{sec:triangulations}
\subsection{Basic definitions}
A \emph{triangulation} of a topological space $X$ is a homeomorphism from $X$ to a simplicial complex. Let us recall that a \emph{simplicial complex} $K$ is specified by a finite set of vertices $V$ and a finite set of simplices $S\subset \mathcal{P}(V)$ (the power set of $V$), such that if
$\sigma \in S$ and $\tau \subset \sigma$ then  $\tau \in S$. The combinatorial data $(V, S)$ is called an {\em abstract simplicial complex}. To each such data, there is an associated topological space, called the {\em geometric realization}. This is constructed inductively on $d \geq 0$, by attaching a $d$-dimensional simplex $\Delta^d$ for each element $\sigma \in S$ of cardinality $d$; see \cite{hat02}. The result is the simplicial complex $K$. In practice, we will not distinguish between $K$ and the data $(S, V)$.

Let $K=(V, S)$ be a simplicial complex. Formally, for a subset $S'\subset S$, its \emph{closure} is $$Cl(S') = \{\tau\in S |\tau \subseteq \sigma\in S'\}$$
The \emph{star} of a simplex $\tau \in S$ is
$$St(\tau) = \{\sigma\in S| \tau\subseteq\sigma\}$$
The \emph{link} of a simplex $\tau \in S$ is
$$Lk(\tau) = \{\sigma\in Cl(St(\tau))| \tau \cap \sigma= \emptyset\}$$

\begin{example} Let $K=\{V,S\}$, where $V=\{1,2,3,4\}$, and
 $$S=\{\{1\}, \{2\}, \{3\}, \{4\}, \{1,2\}, \{1,3\}, \{1,4\}, \{3,4\}, \{1,3,4\} \}.$$ The geometric realization is

\begin{center}
\includegraphics{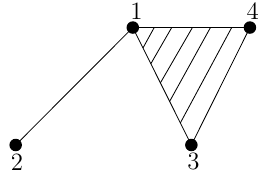}
\end{center}

The link of $\{4\}$ is the edge $\{1,3\}$ (including its vertices, of course). The link of $\{1\}$ is the union of $\{2\}$ and the edge $\{3,4\}$.
\end{example}

Let us mention here that not every CW-complex can be triangulated, but every CW-complex obtained by gluing polyhedral cells in a nice way (by piecewise linear maps) is triangulable \cite{hat02}. See the following example of the torus.

\begin{center}
\begin{picture}(0,0)%
\includegraphics{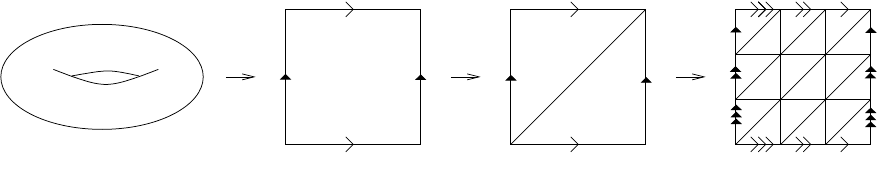}%
\end{picture}%
\setlength{\unitlength}{3158sp}%
\begingroup\makeatletter\ifx\SetFigFont\undefined%
\gdef\SetFigFont#1#2#3#4#5{%
  \reset@font\fontsize{#1}{#2pt}%
  \fontfamily{#3}\fontseries{#4}\fontshape{#5}%
  \selectfont}%
\fi\endgroup%
\begin{picture}(8769,1722)(2093,-2146)
\put(9473,-2082){\makebox(0,0)[lb]{\smash{{\SetFigFont{8}{9.6}{\rmdefault}{\mddefault}{\updefault}{\color[rgb]{0,0,0}simplicial complex}%
}}}}
\put(7426,-2086){\makebox(0,0)[lb]{\smash{{\SetFigFont{8}{9.6}{\rmdefault}{\mddefault}{\updefault}{\color[rgb]{0,0,0}$\Delta$-complex}%
}}}}
\end{picture}%

\end{center}

\subsection{Triangulations of manifolds} \label{sec:manifolds}
In topology, manifolds are considered in different categories, with respect to their transition functions. For example, we have
\begin{itemize}
\item {\em Topological manifolds} if the transition functions are $C^0$;
\item {\em Smooth manifolds} if the transition functions are $C^{\infty}$;
\item {\em PL (piecewise linear) manifolds} if the transition functions are piecewise linear.
\end{itemize}

We say that a triangulation is \emph{combinatorial} if the link of every simplex (or, equivalently, of every vertex) is piecewise-linearly homeomorphic to a sphere. Clearly, every space that admits a combinatorial triangulation is a manifold (in fact, a PL manifold).
\begin{center}
\includegraphics{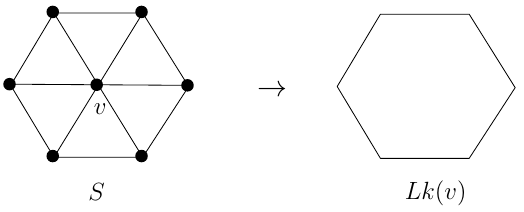}
\end{center}

Conversely, every PL manifold can be shown to admit a combinatorial triangulation.\\

Here are the three main triangulation problems for manifolds.

\subsection*{Question 1} ({Poincar\'e} \cite{poi99})
Does every smooth manifold admit a triangulation?
\subsection*{Answer} (Cairns \cite{cai35}  and Whitehead \cite{whi40}) Yes. Every smooth manifold has an essentially unique PL-structure, and therefore it is triangulable.
\subsection*{Question 2} (Kneser \cite{kne26}) 
Does every topological manifold admit a triangulation?
\subsection* {Answer} This depends on the dimension $n$ of the manifold. \begin{itemize}
\item for $n=0, 1$: Yes, trivially.
\item for $n=2$ (Rad\'o \cite{rad25}) Yes. Every two-dimensional surface has a piecewise linear structure therefore it is triangulable.
\item for $n=3$ (Moise \cite{moi52}) Yes. Every three dimensional manifold is smooth, hence  piecewise linear, and hence triangulable.
\item for $n=4$ (Casson \cite{am90}) No. The Casson invariant can be used to show that Freedman's $E_8$ manifold is not triangulable.
\item for $n\geq5$ (Manolescu \cite{man13}) No. Galewski-Stern \cite{gs80} and Matumoto \cite{mat78} reduced this question to a problem in $3+1$ dimensions. The solution of this reduced problem is given in \cite{man13}, and uses $\pin$-equivariant Seiberg-Witten Floer homology.
\end{itemize}

\subsection*{Question 3} Does every topological manifold admit a PL structure? 
\subsection*{Answer} 
\begin{itemize}
\item for $n\leq3$: Yes, as above.
\item for $n=4$ (Freedman \cite{fre82}) No. Freedman constructed the four-manifold $E_8$ which has no piecewise linear structure.
\item for $n\geq5$ (Kirby-Siebenmann \cite{ks69}) No. For a topological manifold $M$ there exists an obstruction class to having a PL structure. This is called the \emph{Kirby-Siebenmann class} $$\Delta(M)\in H^4(M,\mathbb{Z}/2).$$ 
Vanishing of $\Delta(M)$ is a necessary (and, in dimensions $n\geq5$, sufficient) condition for the existence of a PL structure. In dimensions $n\geq5$ there exist manifolds with $\Delta(M) \neq 0$; e.g. $M=T^{n-4} \times E_8$.\\
\end{itemize}

Note that for dimension $4$ there are more obstructions to the existence of PL (or, equivalently, smooth) structures, apart from $\Delta$. These come from gauge theory; an example is Donaldson's diagonalizability theorem \cite{don83}.

\begin{example} An example of a non-PL triangulation of a manifold can be constructed as follows. Start with a triangulation of a homology sphere $X$ where $\pi
_1(X)\neq 1$. (The Poincar\'e sphere will do.) The suspension $\Sigma X$ is not a manifold, but the double suspension $\Sigma^2 X$ is a topological manifold homeomorphic to a sphere, by the double suspension theorem of Edwards \cite{edw06} , \cite{edw80} and Cannon \cite{can79}. Any triangulation of $\Sigma^2 X$ induced by $X$ is not a combinatorial triangulation. Indeed, the link of any cone point of $\Sigma^2 X$ is $Lk(v)=\Sigma X$, which is not a manifold so is not a PL sphere.\end{example}

Let us mention a few related facts:
\begin{enumerate}
\item[1)]Any manifold of dim $n\neq 4$ is homeomorphic to a CW complex
(Kirby-Siebenmann for dim $n \neq 5$, Quinn for dim $n=5$).
\item[2)]Any manifold of dim $n$ is homotopy equivalent to an $n$-dimensional simplicial complex.\footnote{The fact that the simplicial complex can be taken to be of the same dimension was mistakenly listed as an open problem in the published version of this article.}
\end{enumerate}

Laurence Taylor informed the author of the following proof of 2). Using \cite[Theorem 2C.5]{hat02}, it suffices to check that the manifold is homotopy equivalent to a CW complex of the same dimension. When the dimension $n \neq 4$, this follows from 1) above. When $n=4$ we distinguish several cases:
\begin{enumerate}
\item[(i)] If the manifold $M$ is open then it is smooth, and hence homeomorphic to CW complex;
\item[(ii)] If $M$ is compact with non-trivial boundary, it is homotopy equivalent to its interior, and we can apply Case (i);
\item[(iii)] If $M$ is closed (of any dimension $n$) then it is homotopy equivalent to a simple Poincar\'e complex (roughly, a finite CW complex satisfying Poincar\'e duality). Further, when $n \neq 2$, that simple Poincar\'e complex is homotopy equivalent to a CW complex of dimension $n$; cf. \cite[Theorem 2.2]{wal67}.
\end{enumerate}
This completes the proof.

\smallskip
The following remains unknown:
\begin{enumerate}
\item[1)]Is every $4$-manifold is homeomorphic to a CW complex?
\end{enumerate}

\subsection{The Kirby-Siebenmann obstruction}
Let $M^n$ be a topological manifold of dimension $n\geq 5$. Consider the diagonal $D \subset M \times M $. It can be shown that a small neighborhood $\nu(D)$ of $D$ is an $\R^n$ bundle over $D \cong M$. This is called \emph{topological tangent bundle} of $M$.

One can construct some infinite dimensional topological groups as follows: \\

\noindent $ TOP = \underrightarrow{\lim} \ TOP(n),$ where $TOP(n)$ consists of the homeomorphisms of $\R^n$ fixing $0$;

\noindent $ PL = \underrightarrow{\lim} \ PL(n) \subset TOP,$ where $PL(n)$ consists of the PL-homeom. of $\R^n$ fixing $0$.\\
 
 There is a fibration
 \begin{center}
 \begin{tikzcd}
 K(\mathbb{Z}/2,3) = TOP/PL\arrow{r}  & BPL \arrow{d}{\Psi}  \\
 {} & BTOP 
 \end{tikzcd} 
  \end{center}

Recall that in Question 3 above, we mentioned that the obstruction to the existence of a PL structure on $M$ is the Kirby-Siebenmann class $\Delta(M)\in H^4(M,\mathbb{Z}/2)$. For general $M$, the class $\Delta(M)$ is defined as the obstruction to lifting the map in the diagram below:

\begin{center}

\begin{tikzcd}
 {} & BPL \arrow{d}{\Psi}  \\
M \arrow [dashed]{ur}\arrow{r}{\Phi} & BTOP
\end{tikzcd}

\end{center}

\vspace{0.3cm} If $M^n$ has a triangulation $K$ (not necessarily PL), we can give a more concrete definition of the Kirby-Siebenmann class $\Delta(M)$, as follows. For simplicity, let us assume that $M$ is orientable. Let
\begin{equation}
\label{eq:cK}
 c(K)=\sum_{\sigma\in K^{n-4}} [Lk(\sigma)]\sigma \in H_{n-4}(M,\Theta_3 ^H) \cong H^{4}(M,\Theta_3 ^H)
 \end{equation}

\begin{center}
\includegraphics{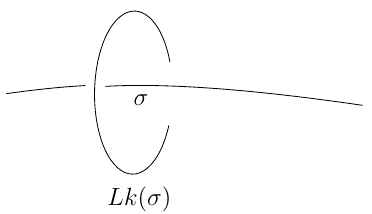}
\end{center}
If $M$ is non-orientable, we can use Poincar\'e duality with local coefficients and still obtain an element $c(K) \in H^4(M; \Theta^3_H)$.

Here, $\Theta_3 ^H$ represents the {\em homology cobordism group} of homology three-spheres:
$$\Theta_3 ^H= \langle Y^3 \text{ oriented } \mathbb{Z}HS^3 \rangle/\sim$$
where the equivalence relation is
$Y_0\sim Y_1 \iff \exists W^4 \text{ (PL or, equivalently, smooth) }$ such that $\partial(W)= Y_0 \cup (- Y_1)$ and $H_*(W,Y_i;\mathbb{Z})=0$. Addition in $\Theta_3 ^H$ is connected sum and the identity element is $[S^3]=0$. 
\begin{center}
\includegraphics{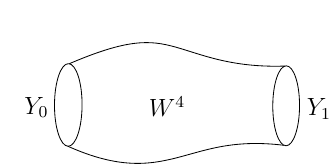}
\end{center}

\bigskip	  
The structure of the Abelian group $\Theta_3 ^H$ is not fully understood. Fintushel-Stern \cite{fs85} \cite{fs90} and Furuta \cite{fur90} proved that $\Theta_3 ^H$ is infinitely generated. Fr{\o}yshov \cite{fro10} proved that $\Theta_3 ^H$ has a $\mathbb{Z}$-summand. However, the following questions remain open.
	
\begin{itemize}
		\item Does $\Theta_3 ^H $ have a torsion? 
		\item Does $\Theta_3 ^H $ have a $\mathbb{Z}^\infty$-summand? 
\end{itemize}

	The analogous homology cobordism groups $\Theta_n ^H$ are all trivial for $n\neq3$, by the work of Kervaire \cite{ker69}. 
	
	To see that the homology cobordism group in dimension three is nontrivial ($\Theta_3 ^H\neq 0$), one can use the \emph{Rokhlin homomorphism}, defined using the \emph{Rokhlin invariant}  $$\mu = \Theta_3 ^H \rightarrow \mathbb{Z}/2, \ \ \  \mu(Y)= \mu(Y, \s) = \frac{\sigma(W)}{8} \pmod{2}. $$ 
	Here, $(W,\mathfrak{t})$ is an arbitrary compact, smooth, $spin(4)$ manifold with $\partial(W,\mathfrak{t})= (Y,\mathfrak{s})$ and $\sigma(W)$ denotes the signature of $W$. When $Y$ is an integral homology sphere, there is a unique $spin$ structure $\mathfrak{s}$ on $Y$. 
	
	For example, we have $\mu(S^3)=0$ and $\mu(Y)=1$, where $Y$ is the Poincar\'e homology sphere. We get that $\Theta_3 ^H\neq 0$.
	
	Coming back to the element $c(K)$ from \eqref{eq:cK}, consider the following short exact sequence and the associated long exact sequence on cohomology.
	
	\begin{equation}\label{1}
	0\longrightarrow Ker\mu\longrightarrow \Theta_3 ^H \stackrel{\mu}{\longrightarrow} \mathbb{Z}/2\longrightarrow 0
	\end{equation}

\begin{eqnarray}
\label{eq:les}
		{} \dots \rightarrow & H^4(M;\Theta_3 ^H)  & \stackrel{\mu}{\longrightarrow}  H^4(M;\mathbb{Z}/2) \stackrel{\delta}{\longrightarrow} H^5(M;Ker\mu)\longrightarrow \dots \\
		{} & c(K) & \stackrel{\mu}{\longmapsto } \Delta(M) \notag
\end{eqnarray}
It can be shown that the image of $c(K)$ under $\mu$ is exactly the Kirby-Siebenmann class. Note that $c(K)$ is zero if the triangulation $K$ is combinatorial. By contrast, $\mu(c(K)) = \Delta(M)$ is zero if and only if $M$ admits some combinatorial triangulation (possibly different from $K$).

\subsection{Triangulability of manifolds}
If the manifold $M$ has a triangulation, from the fact that \eqref{eq:les} is a long exact sequence we see that
$$\delta(\Delta(M))=0 \in  H^5(M;Ker\mu).$$ 
By Galewski-Stern \cite{gs80} and Matumoto \cite{mat78}, the converse is also true. More precisely,  for dim $M = n\geq 5$,
$$M^n \text{ is triangulable} \iff\delta(\Delta(M))=0. $$
Furthermore, they show that
$$(\ref{1}) \text{ does not split } \iff \forall n\geq 5, \ \exists M^n{, } \ \delta(\Delta(M))\neq 0.$$

\begin{theorem}[\cite{man13}]
\label{thm:man}
The short exact sequence (\ref{1}) does not split.
\end{theorem}

Theorem~\ref{thm:tc} is a consequence of Theorem~\ref{thm:man}, combined with the work of Galewski-Stern and Matumoto. 

To construct explicit examples of non-triangulable manifolds, consider the short exact sequence
$$0\longrightarrow \mathbb{Z}/2 \longrightarrow \mathbb{Z}/4 \longrightarrow \mathbb{Z}/2\longrightarrow 0.$$
The Bockstein (connecting) homomorphism associated to this sequence is the  first Steenrod square on cohomology,
$$H^k(M;\mathbb{Z}/2) \stackrel{Sq^1}{\longrightarrow} H^{k+1}(M;\mathbb{Z}/2)$$  

\begin{exercise}If $M$ is a manifold of dimension $\geq 5$ and $Sq^1\Delta (M)\neq 0$, show that $\delta(\Delta(M))\neq 0$. {\em Hint:} Make use of the fact that (\ref{1}) does not split.
\end{exercise}

Thus, it suffices to find a $5$-dimensional example $M^5$ such that $Sq^1(\Delta(M))\neq 0$. It will then follow that $M$, and hence $M \times T^{n-5}$ for $n \geq 5$, are non-triangulable.

\begin{example} (Kronheimer) Let $X=*(\mathbb{
CP}^2\# \overline{\mathbb{CP}^2})$ be a simply connected topological $4$-manifold with intersection form $$\begin{pmatrix} 
 1 & 0 \\
 0 & -1 
\end{pmatrix} \sim - \begin{pmatrix} 
  1 & 0 \\
   0 & -1 
\end{pmatrix}$$ and  $\Delta(M) \neq 0$. Such an $M$ exists by Freedman's work \cite{fre82}. Moreover, Freedman's theory also shows the existence of an  orientation reversing homeomorphism $f:X\to X $. Let $M^5$ be a mapping torus of $f$
$$M={(X\times I)}/{(x,0)\sim(f(x),0)}.$$
We have
$\Delta(X)=1 \in H^4(X;\mathbb{Z}/2)=\mathbb{Z}/2$, and therefore
$Sq^1\Delta(M)=\Delta(M)\cup w_1(M)\neq 0.$
\end{example}

In fact, all non-triangulable manifolds of dim $n=5$ are non-orientable. In dim $n\geq 6$ there also exist orientable examples. The simplest such example is $P^6$, the circle bundle over the manifold $M^5$ from the example above, associated to the oriented double cover of $M^5$.

Let us finish this section with a Venn diagram showing the different kinds of manifolds.\\

\begin{center}
\includegraphics{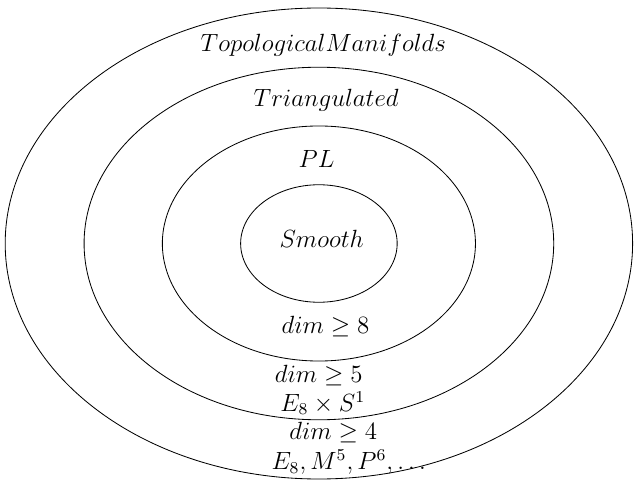}
\end{center}

\medskip
In the lower part of the diagram we indicated the lowest dimension for which the respective set difference is non-empty, and some examples of manifolds with the required properties.

\section{Seiberg-Witten theory}
\label{sec:sw}

In Sections~\ref{sec:sw} through \ref{sec:duality}, we will sketch the proof of Theorem~\ref{thm:man}, following \cite{man13}.

\subsection{Strategy of proof} \label{sec:strategy}
We want to show that the exact sequence 
$$	0\longrightarrow Ker\mu\longrightarrow \Theta_3 ^H \stackrel{\mu}{\longrightarrow} \mathbb{Z}/2\longrightarrow 0 $$
does not split. This is equivalent to proving that there is no $Y \in\mathbb{Z}HS^3$ such that $2[Y]=0$ in $\Theta_3 ^H$ and $\mu(Y)=1$. Obviously $2[Y]=0$ means $Y\sim -Y$ or, equivalently, $Y\#Y\sim S^3$.

The idea is to construct a map $$\mathbb{Z}HS^3\stackrel{\beta}{\longrightarrow} \mathbb{Z}$$ with certain properties. Specifically, for a homology sphere $Y$, we define an invariant $\swfh_*^{\pin}(Y)$, called the \emph{$\pin$-Equivariant  Seiberg-Witten Floer Homology}, and then set
$$\beta(Y)=\frac{1}{2}\bigl( (\text{min.deg.of middle } v\text{-tower in } \swfh_*^{\pin}(Y))-1\bigr ).$$
We will explain what the middle $v$-tower is later. For now, let us state the properties of $\beta$ that are of importance to us:
\begin{enumerate}
\item[1)]  $\beta$ descends to $\Theta_3 ^H$, i.e. $Y_0\sim Y_1 \implies \beta(Y_0)=\beta(Y_1)$
\item[2)]$ \mu(Y)= \beta(Y) \ (\text{mod }2)$, where $\mu$ is the Rokhlin invariant;
\item[3)]$ \beta(-Y)=-\beta(Y)$.
\end{enumerate}

It should be noted that $\beta$ is not a homomorphism. For example, for the Brieskorn sphere $Y=\Sigma(2,3,11)$, one can compute that $\beta(Y\#Y)=2$ while $\beta(Y)=0$. 

Once we construct $\beta$ with the properties above, the argument below will give the result of  Theorem~\ref{thm:man}.
\begin{eqnarray*}
& \text{If } \exists Y\in\mathbb{Z}HS^3\text{ such that }Y\sim -Y \\
& \implies \beta(Y)= \beta(-Y)=-\beta(Y)\\
& \implies \beta(Y)=0 \\
& \implies \mu(Y)=0. \text{ } 
\end{eqnarray*}
This contradicts that $\mu(Y)\neq1$.

The construction of $\beta$ is inspired by the following, previously defined invariants:
\begin{itemize}
\item The Casson invariant $\lambda(Y)$ from \cite{am90}, which satisfies properties 2) and 3) but not 1).
\item The Froysh\o v invariant $\delta(Y)=\frac{1}{2}( \text{min.deg.of } U\text{-tower in } \swfh_*^{S^1}(Y))$, and the Ozsv\'ath-Szab\'o correction term $d(Y)$ (conjecturally, $\delta=\frac{d}{2}$). See \cite{fro10}, \cite{os03}. These invariants satisfy properties 1) and 3) but not 2).
\end{itemize}

In the definition of the Fr{\o}yshov invariant, $\swfh_*^{S^1}(Y)$ is the {\em $S^1$-Equivariant Seiberg-Witten Floer Homology} of $Y$. Variants of this theory were constructed by Marcolli-Wang in \cite{mw01}, the author in \cite{man03}, and Fr{\o}yshov in \cite{fro10}, all for the case of manifolds with $b_1(Y)=0$ (rational homology spheres), as well as by Kronheimer-Mrowka \cite{km07} for all three-manifolds.

The invariant $\beta$ is the analogue of $\delta$, but using $\pin$-equivariant, instead of $S^1$-equivariant, Seiberg-Witten Floer homology. The $\pin$-equivariant theory was first defined by the author in \cite{man13}, for rational homology three-spheres, and this is the construction that we will review in these notes. Since then, a different construction of the $\pin$-equivariant theory was given by Lin in \cite{lin15}, and that applies to all three-manifolds. Lin's construction provides an alternative proof of Theorem~\ref{thm:tc}.

\subsection{The Seiberg-Witten equations}
The constructions below work (with minor modifications) for all rational homology $3$-spheres. However, for simplicity, and since this is what we need for Theorem~\ref{thm:man}, we will only discuss the case when $Y$ is an integral homology $3$-sphere.

Let $Y$ be such a homology sphere, and pick a Riemannian metric $g$ on $Y$. There is a unique $Spin^c$ structure on $(Y,g)$, denoted $\mathfrak{s}$. Specifically, $\mathfrak{s}$ consists of a rank $2$ Hermitian bundle $\mathbb{S}$ on $Y$, together with a Clifford multiplication map 
$$\rho : TY\stackrel{\cong}{\longrightarrow}\mathfrak{su}(\mathbb{S})\subseteq End(\mathbb{S}).$$ 

Here, $\mathfrak{su}(\mathbb{S})$ denotes traceless ($trA=0$), skew-adjoint($A+A^*=0)$ endomorphisms of $S$. Similarly, below, $\mathfrak{sl}(\mathbb{S})$ will denote traceless ($trA=0$) endomorphisms of $\mathbb{S}$. After complexification, and using the duality $TY\cong T^*Y$, the Clifford multiplication extends to an isomorphism
\begin{equation}
\label{eq:rho}
\rho=T^*Y\otimes\mathbb{C}\stackrel{\cong}{\longrightarrow}\mathfrak{sl}(\mathbb{S})\subseteq End(\mathbb{S})
\end{equation}

Explicitly, we can construct the $Spin^c$ structure as follows. We let $\mathbb{S}$ be a trivial bundle $\mathbb{S}=\underline{\mathbb{C}}^2$. Trivialize 
$$TY= \langle e_1, e_2, e_3 \rangle$$ and define $\rho$ by
\begin{center}
$\rho(e_1)=\begin{pmatrix} 
 i & 0 \\
 0 & -i 
\end{pmatrix}$, \hspace{.5cm} $\rho(e_2)=\begin{pmatrix} 
 0 & -1 \\
 1 & 0 
\end{pmatrix}$, \hspace{.5cm} $\rho(e_3)=\begin{pmatrix} 
 0 & i \\
 i & 0 
\end{pmatrix}.$
\end{center}

There is an associated {\em Dirac operator}
$${\partiali}:\Gamma(\mathbb{S})\mapsto\Gamma(\mathbb{S}), \hspace{1cm} \partiali(\phi)=\sum_{1}^{3}\rho(e_i)\partial_i(\phi) \hspace{1cm} (\partial_i\leftrightarrow e_i)$$

Further, we have an identification
\begin{eqnarray*}
\{Spin^c \text{ connections on } \mathbb{S} \} &\stackrel{\cong}{\leftrightarrow}& \Omega^1(Y,i\mathbb{R}) \\
A&\leftrightarrow& a\\
A & = & A_{0}+a,
\end{eqnarray*}
where $A_{0}$ is the trivial connection on $\mathbb{S}$.

For $\phi\in\Gamma(\mathbb{S})$, consider the endomorphism $(\phi\otimes\phi^*)_{\circ}\in \mathfrak{sl}(S)$, the traceless part of $(\phi\otimes\phi^*)$. Using \eqref{eq:rho}, we get a form
$$\rho^{-1}((\phi\otimes\phi^*)_{\circ})\in \Omega^1(Y;\mathbb{C}).$$

 Consider the configuration space $\mathcal{C}(Y,\mathfrak{s})=\Omega^1(Y;i\mathbb{R})\oplus \Gamma(\mathbb{S})$. For a pair $(a,\phi)\in \mathcal{C}(Y,\mathfrak{s})$, the {\em Seiberg-Witten Equations} are 

\begin{align}
 \widetilde{SW}(a,\phi) = \left\{\begin{array}{c l} 
*da -2\rho^{-1} ((\phi\otimes\phi^*)_{\circ})&=0 \\
\partiali\phi+\rho(a)\phi&=0 \end{array}\right.
\end{align}

In a short way we write this as  $$\widetilde{SW}(a,\phi) = 0.$$

The Seiberg-Witten equations are invariant under the action of the gauge group $\mathcal{G}= C^\infty(Y,S^1)=\{u:Y\to S^1\}$,
$$ u \cdot (a,\phi)=(a-u^{-1}du,u\phi) $$

Since $Y$ is a homology sphere, any $u: Y \to S^1$ can be written as $u=e^\xi$ for some $\xi:Y\to i\mathbb{R}$. The action becomes
$$ e^{\xi} \cdot (a, \phi) = (a-d\xi, e^\xi\phi).$$

We can think of the Seiberg-Witten equations as the gradient flow equations for the {\em Chern-Simons-Dirac functional}
$$CSD(a,\phi)=\frac{1}{2}\big( \int_{Y}{\langle \phi,(\partiali\phi+\rho(a)\phi) \rangle dvol}-\int_{Y} a\wedge da\big)$$
We have $\widetilde{SW}=\nabla{(CSD)}$.

\subsection{The Seiberg-Witten equations in Coulomb gauge} 
Define the (global) \emph{Coulomb slice} $$V:=ker(d^*)\oplus \Gamma(S)\subset\mathcal{C}(Y,\mathfrak{s}).$$
In other words, we restrict the configuration space $\mathcal{C}(Y,\mathfrak{s})=i\Omega^1(Y;i\mathbb{R})\oplus \Gamma(\mathbb{S})$ by adding the condition $d^*a=0$. We can view $V$ as as the quotient of the configuration space by the \emph{normalized gauge group action} of $\mathcal{G}_0\subset \mathcal{G}$. Here, 
$$\mathcal{G}_0 = \{ u : Y \to S^1 \mid u=e^\xi, \ \int_Y \xi = 0 \}.$$

Since $Y$ is a homology sphere, we have a Hodge decomposition $$\Omega^1(Y)=ker(d)\oplus ker(d^*).$$
At $(a, \phi) \in V$, let $\pi_V$ denote the linear projection from the tangent space $T_{(a, \phi)} \Conf(Y, \s)$ onto $V$, with kernel the tangents to the $\mathcal{G}_0$-orbit. Note that $\pi_V$ is {\em not} an $L^2$-orthogonal projection! 

Let
$$SW:V\to V, \hspace{2cm}SW:=\pi_{V}\circ \widetilde{SW}.$$

Using $e^{i\theta}:(a,\phi)\mapsto(a,e^{i\theta }\phi)$ we find a bijection:
$$\{\text{Flow lines of }\widetilde{SW}\}\big/ \mathcal{G}\stackrel{1:1}{\longleftrightarrow}\{\text{Flow lines of }SW \}\big/ S^1. $$ 

Furthermore, $V$ has a metric $\tilde{g}$ induced by  
$$\langle v,w \rangle_{\tilde{g}}= \langle \pi^{elc}(a),\pi^{elc}(b) \rangle_{L^2},$$
where $\pi^{elc}$ is the $L^2$-orthogonal projection from $T_{(a, \phi)} \Conf(Y, \s)$ with kernel the tangent to the $\mathcal{G}_0$-orbit. The image of $\pi^{elc}$ is an {\em enlarged local Coulomb slice} $K^{elc}$, the orthogonal complement to the $\mathcal{G}_0$-orbit. See the figure below.

\vspace{.2cm}
\begin{center}
\begin{picture}(0,0)%
\includegraphics{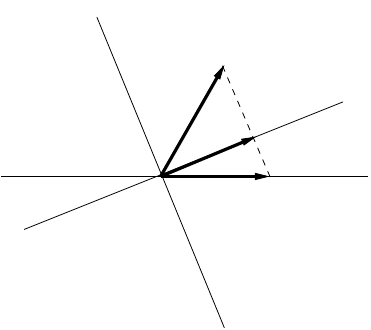}%
\end{picture}%
\setlength{\unitlength}{3355sp}%
\begingroup\makeatletter\ifx\SetFigFont\undefined%
\gdef\SetFigFont#1#2#3#4#5{%
  \reset@font\fontsize{#1}{#2pt}%
  \fontfamily{#3}\fontseries{#4}\fontshape{#5}%
  \selectfont}%
\fi\endgroup%
\begin{picture}(3474,3096)(1489,-4498)
\put(4881,-3238){\makebox(0,0)[lb]{\smash{{\SetFigFont{10}{12.0}{\rmdefault}{\mddefault}{\updefault}{\color[rgb]{0,0,0}$V$}%
}}}}
\put(2551,-1561){\makebox(0,0)[lb]{\smash{{\SetFigFont{10}{12.0}{\rmdefault}{\mddefault}{\updefault}{\color[rgb]{0,0,0}$T\mathcal{G}_0$}%
}}}}
\put(4276,-2386){\makebox(0,0)[lb]{\smash{{\SetFigFont{10}{12.0}{\rmdefault}{\mddefault}{\updefault}{\color[rgb]{0,0,0}$K^{elc}$}%
}}}}
\put(3301,-2236){\makebox(0,0)[lb]{\smash{{\SetFigFont{10}{12.0}{\rmdefault}{\mddefault}{\updefault}{\color[rgb]{0,0,0}$v$}%
}}}}
\put(3931,-2807){\makebox(0,0)[lb]{\smash{{\SetFigFont{10}{12.0}{\rmdefault}{\mddefault}{\updefault}{\color[rgb]{0,0,0}$\pi^{elc}(v)$}%
}}}}
\put(3762,-3259){\makebox(0,0)[lb]{\smash{{\SetFigFont{10}{12.0}{\rmdefault}{\mddefault}{\updefault}{\color[rgb]{0,0,0}$\pi_V(v)$}%
}}}}
\end{picture}%

\end{center}
\medskip

On the Coulomb slice $V$, the $SW$ equation can be written as the sum of a linear part and another part
$$SW = l+c,$$ 
where $l,c:V\to V$ are given by 
$$\begin{array}{c l}
l(a,\phi) =& (*da,\partiali\phi)\\
c(a,\phi) =& \pi_V\circ(-2\rho^{-1}(\phi\otimes\phi^*)_\circ, \rho(a)\phi)
\end{array}$$

Let $V_{(k)}$ be the $L_k^2$-completion of $V$ for a fixed number $k\gg 0$. We will take $k>5$. Then, the map $l:V_{(k)}\to V_{(k-1)}$ is a linear, self-adjoint, Fredholm operator, and $c:V_{(k)}\to V_{(k-1)}$ is a compact operator.

The following is the standard compactness theorem for Seiberg-Witten equations, adapted to Coulomb gauge.
\begin{theorem}
\label{thm:cpt}
 Fix $k>5$. There exists some $R > 0$ such that all the critical points and flow lines between critical points of $SW$ are contained inside the ball $\mathcal{B}(R)\subset V_{(k)}$. 
\end{theorem}

\subsection{Finite dimensional approximation}
Seiberg-Witten Floer homology is meant to be Morse homology for the functional $SW$ on $V$. However, instead of finding a generic perturbation to achieve transversality, it is more convenient to do finite dimensional approximation. In the finite  dimensional case, we can simply use singular homology instead of Morse homology.

Our finite dimensional approximation is inspired by Furuta's $4$-dimensional case. In our setting, $V$ is an infinite dimensional space, and as a finite dimensional approximation of $V$ we consider 
$$V_\lambda^\mu=\oplus \ \bigl (\text{eigenspaces of } l \text{ with eigenvalues in } (\lambda,\mu) \ \bigr ), \ \ \ \lambda \ll 0 \ll \mu.$$
 
We replace $SW=l+c:V\to V$  by  
$$l+p_\lambda^\mu c: V_\lambda^\mu\to V_\lambda^\mu,$$ 
where $p_\lambda^\mu$ is the $L^2$ projection onto $V_\lambda^\mu$. Then,
 $$SW_\lambda^\mu=l+p_\lambda^\mu c$$ is a vector field on $V_\lambda^\mu$.

The following is a compactness theorem in the finite dimensional approximations.

\begin{theorem}[cf. \cite{man03}]
There exists $R>0$ such that for all $\mu\gg0\gg \lambda$ all critical points of $SW_\lambda^\mu$ in the ball $\mathcal{B}(2R)$, and all flow lines between them that lie in $\mathcal{B}(2R)$, actually  stay in the smaller ball $\mathcal{B}(R)$.\end{theorem}

\begin{center}
\includegraphics{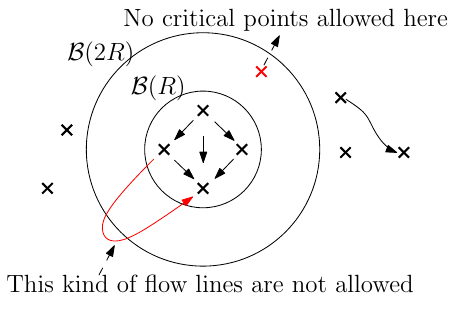}
\end{center}

The idea of the proof is to use that in $\mathcal{B}(2R)$, $l+p_\lambda^\mu c$ converges to $l+c$ uniformly. We then apply Theorem~\ref{thm:cpt}.

\subsection{The Conley index}
\label{sec:Conley}
In finite dimensions, on a compact manifold Morse homology is just the ordinary homology. In our setting we work on the non-compact space $\mathcal{B}(2R) \subset V_{\lambda}^{\mu}$. On a non-compact manifold, Morse homology is the homology of a space called the {\em Conley index}.

The Conley index (cf. \cite{con78}) is defined for any isolated invariant set $S$ of a flow $\{\phi_t\}$ on an $m$-dimensional manifold $M$. 

\begin{definition} For a subset $A \subseteq M$ we define 
$$Inv\; A=\{x\in M \mid \phi_t(x)\in A, \forall t\in \mathbb{R} \}$$
\end{definition}

\begin{definition}
A compact subset $S\subset M$ is called an \emph{isolated invariant set} if there exists $A$ a compact neighborhood of $S$ such that $S=Inv \; A \subseteq Int \; A$.
\end{definition}

\begin{definition}
For an isolated invariant set $S$, the \emph{Conley index} is defined as $I(S)=N/L$  where $L\subseteq N\subseteq M$, with both $L$ and $N$ compact, satisfy
\begin{enumerate}
\item[1)] $Inv(N-L)=S\subseteq Int(N-L)$;
\item[2)] $\forall x\in N$, if $\exists t>0 \text{ such that } \phi_t(x)\notin N$,  then $\exists \tau \in [0, t)$  with $\phi_{\tau}(x) \in L$;
\item[3)] $x\in L, t>0, \phi_{[0,t]}(x) \subset N \implies \phi_{[0,t]}(x)\subset L$.
\end{enumerate}
\end{definition}

\medskip
\begin{center}
\includegraphics{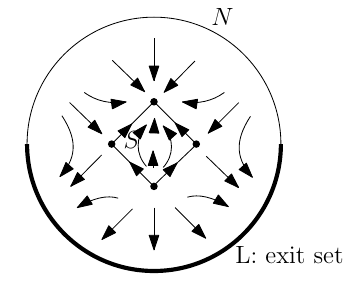}
\end{center}
\medskip

In our case we take $A=\mathcal{B}(2R)$, so that $S=Inv \; A$ (the union of critical points and flow lines inside the ball). Then, $N$ can be taken to be a manifold with boundary and $L\subset\partial N$ a codimension $0$ submanifold of the boundary (so that $L$ has its own boundary).

It can be shown that, if the flow lines satisfy the Morse-Smale condition, then Morse homology is isomorphic to the reduced singular homology of $I(S)$.

\begin{example}
As an example of Conley index, let $S=\{x\}$ be an index $k$ Morse critical point. We can choose $N=D^k\times D^{n-k}$ and $L=\partial D^k\times D^{n-k}$:
\begin{center}
\includegraphics{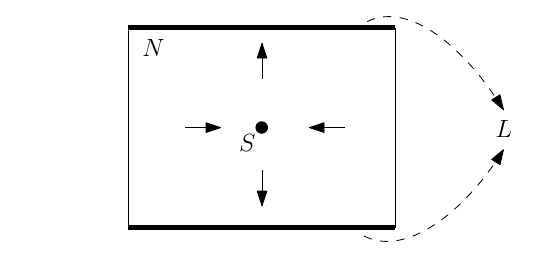}
\end{center}
The Conley index is then a sphere, $N\big/L\simeq S^k$.
\end{example}

\subsection{Seiberg-Witten Floer Homology}
Following \cite{man03}, we define the $S^1$-equivariant Seiberg-Witten Floer homology of a homology sphere $Y$ to be
$$\swfh_*^{S^1}(Y)=\widetilde{H}_{*+shift}^{S^1}(I_{\lambda}^\mu) \text{ for } \mu\gg0\gg \lambda.$$ 

Here, $I_{\lambda}^\mu$ is Conley index for $S^\mu \subset V_\lambda^\mu $, and the grading is shifted by some amount that depends on $\lambda$, $\mu$, and the metric. (We will specify the precise shift later.) Note that everything in the construction is $S^1$-equivariant, with respect to the $S^1$ action by constant gauge transformation. This allows us to apply Borel homology in the formula above.

In fact, everything is $\pin$-equivariant, because in our situation the $Spin^c$-structure is $Spin$. Here, 
$$\pin= S^1\cup jS^1\subset \mathbb{C}\cup j\mathbb{C}=\mathbb{H}, \ \ j^2=-1 , \ ij=-ji,$$ 
and $j$ acts on $V$ as an additional symmetry 
$$j:(a,\phi)\to (-a,\phi j)$$
$$\Gamma(S)=\{U\mapsto\mathbb{C}^2\}$$
$$\begin{pmatrix}
v_1\\
v_2
\end{pmatrix} \to \begin{pmatrix}
-\bar{v}_2\\
-\bar{v}_1
\end{pmatrix}$$

The Seiberg-Witten equations are $\pin$-equivariant. Thus, following \cite{man13}, we can give the definition of $\pin$-equivariant Seiberg-Witten Floer homology. It is helpful to use coefficients in the field $\mathbb{F}=\mathbb{Z}/2$.
\begin{definition}
$$\swfh_*^{\pin}(Y; \mathbb{F})=\widetilde{H}_{*\,+shift }^{\pin}(I_{\lambda}^\mu; \mathbb{F}) \ \text{ for } \mu\gg0\gg \lambda.$$ 
\end{definition}

\subsection{Invariance}
Let us prove that the Floer homologies we have just defined, $\swfh_*^{S^1}(Y)$ and $\swfh_*^{\pin}(Y)$, are invariants of $Y$. In the process, we will also identify the grading shift in their definitions.

Remember we have a map $SW:V\to V$, where $V=kerd^*\oplus\Gamma(S)$ and $V_\lambda^\mu$ $(\lambda\ll0\ll\mu)$  is a finite dimensional approximation of $V$. The the $SW$-equations can be decomposed as $SW=l+c$, with approximations
$$SW_\lambda^\mu=l+p_\lambda ^\mu c:V_\lambda^\mu\to V_\lambda^\mu.$$ 

The flow equation is 
$$\dot{x}=-SW_\lambda^\mu(x(t)).$$

Let us investigate how the Conley index $I^{\mu}_{\lambda}$ changes under varying $\mu$ and $\lambda$. If we change $\mu \leadsto \mu'>\mu\gg 0 $, we have a decomposition
$$
\begin{array}{cccccc}
V_{\lambda}^{\mu'} &=& V_{\lambda}^{\mu}&\oplus & V_{\mu}^{\mu'}& \\
\vdots & &\vdots& &\vdots & \\
l+p_\lambda ^{\mu'} c&\leadsto & l+p_\lambda ^{\mu} c & \oplus & l+p_\mu ^{\mu'} c& (\simeq l \text{ linear flow}) 
\end{array}
$$

The Conley index invariant under deformations i.e, if we have a family of flows $\varphi(s)$ where $s \in [0,1]$, such that  
$$S{(s)}  =  Inv(\mathcal{B}(R) \text{ in } \! \varphi(s)) \subset Int \; \mathcal{B}(R), \text{ where } s\in[0,1],$$
then $I(S(0)) \simeq  I(S(1)).$ 

\medskip

\begin{center}
\includegraphics{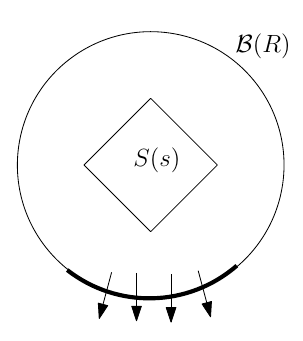}
\end{center}

\medskip

In our case, we let $\varphi(0)$ be the flow of $l + p^{\mu'}_{\lambda}$ and deform it into $\varphi(1)$, the direct sum of the flow of $l+p^{\mu}_{\lambda}$ and the linear flow $l$ on $V^{\mu'}_{\mu}$. We get
$$I_\lambda^{\mu'} =  I(S(0)) = I(S(1)) = I_\lambda^\mu \wedge I_{\mu}^{\mu'}(l).$$ 
Here, $ I_{\mu}^{\mu'}(l)$ is the Conley index for the  linear flow  $\dot{x}=-l(x)$ on $V_\mu^{\mu'}$. Since the restriction of $l$ to that subspace has only positive eigenvalues, we see that
$$ I_{\mu}^{\mu'}(l)= S^{(\text{Morse index})}=S^0.$$

We obtain that
$$  I_\lambda^{\mu'}=I_\lambda^\mu \text{ when } \mu,\mu^{,}\gg 0. $$

On the other hand, by a similar argument, when we vary the cut-off $\lambda$ for negative eigenvalues, the Conley index changes by the formula

$$ \begin{array}{ccccc}
 I_{\lambda^{'}}^{\mu} & = & I_\lambda^\mu &\wedge & \underbrace{ ( {V^\lambda_{\lambda^{'}}} )^{+}}.\\ 
 & & & & \text{sphere}
\end{array} $$

We conclude that:
$$\widetilde{H}_{*\,+\,dim \, V_\lambda^0 }^{S^1}(I_{\lambda}^\mu)$$ 
is independent of $\lambda$ and $\mu$, provided $\mu \gg 0 \gg \lambda$. The same is true for the $\pin$-equivariant homology. 

This suggests including a degree shift $dim \, V_{\lambda}^0$ in the definitions. However, we still have to investigate the dependence on the Riemannian metric $g$.

Fix $\mu, \, \lambda$ such that they are not eigenvalues of $l$. Then, varying $g$ does not change the Conley index.

\begin{center}
 \includegraphics{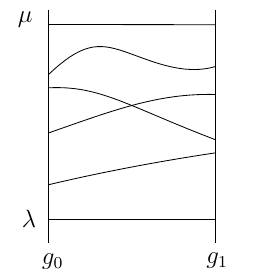}
\end{center}

Note that $V_\lambda^\mu$ has the same dimension as we vary $g$. The problem is that as we vary $g$, the dimension of $V_\lambda^0$ may change. The change is governed by the spectral flow of the linear operator $l$, that is, the signed count of eigenvalues that cross $0$ as we vary the metric. In the picture below, this count is $2-1=1$.
 
 \begin{center}
  \includegraphics{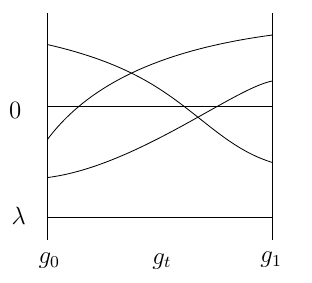}
 \end{center}
 
 Notice that in the linear part  $l=(*d,\Di)$ of the equation, $*d$ does not have zero eigenvalues,  since $H^1=0$. However, $\Di$ has spectral flow. Pick a spin four-manifold $W$ with boundary $(Y, g)$, and add a cylindrical end:
 
 \begin{center}
   \includegraphics{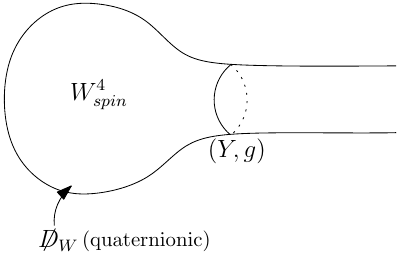}
  \end{center}
 
Then, the spectral flow of $\Di$ is given by the formula 
 \[ \begin{array}{ccl}
 SF(\Di) & = & n(Y,g_0)-n(Y,g_1)\\
  & = & 2\, ind\,(\Di) \text{ on } Y\times[0,1] 
 \end{array} \]
where $$n(Y,g)= -2\,(ind_{\mathbb{C}}(\Di_W)\, +\, \frac{\sigma(W)}{8}) \in 2\mathbb{Z}$$ 
and
 $$n(Y,g)\equiv 2\mu \pmod{4}.$$
 Here, $\mu$ is the Rokhlin invariant.
 
Hence, $$\swfh_*^{S^1}(Y):=\widetilde{H}_{*\,+\,dim \, V_\lambda^0\,-\,n(Y,g) }^{S^1}(I_{\lambda}^\mu)$$ is an invariant of $Y$. The same goes for
$$\swfh_*^{\pin}(Y; \mathbb{F}):=\widetilde{H}_{*\,+\,dim \, V_\lambda^0\,-\,n(Y,g) }^{\pin}(I_{\lambda}^\mu ;\mathbb{F}).$$

One can also do similar constructions with other generalized homology theories, such as $K$-theory $\widetilde{K}_*$ and equivariant $K$-theory $\widetilde{K}_*^{\pin}$, or to the Borel homology $\widetilde{H}_*^G$ for any subgroup $G \subset \pin$, etc.

\subsection{The Seiberg-Witten Floer stable homotopy type}
In fact, the most general invariant produced by the above construction is a $\pin$-equivariant stable homotopy type, denoted $\swf(Y)$. We will refer to it as an equivariant suspension spectrum (although our definition is slightly different from the one of equivariant spectra in the algebraic topology literature).

\begin{definition}
Non-equivariantly, a \emph{suspention spectrum} is a pair $(X,n)$ consisting of a pointed topological space $X$ and $n\in \mathbb{R}.$ We think of $(X, n)$ as the formal de-suspension of $X$, $n$ times: 
$$(X,n)=\Sigma^{-n}X.$$
\end{definition}

Recall that the $n$th suspension of a space $X$ is 
$$\Sigma^{n}\,X =  S^n\wedge X.$$
For example, $\Sigma^{n}\,S^k =  S^{n+k}$. In the world of suspension spectra, we can talk about the $(-n)$-dimensional sphere:
$$(S^0, n) = \Sigma^{-n}S^0 =  S^{-n}.$$

We have a formal suspension $\Sigma$ given by $$\Sigma(X,n) = (X,n+1).$$ 
For any integer $m$, the following identification holds:
$$ \Sigma^m(X,n) \sim (X,n-m) \sim  (\Sigma^mX,n).$$

Let $[X, Y]$ denote the set of homotopy classes of (pointed) maps from $X$ to $Y$. One can define a category with objects and morphisms\\

 $Obj= (X,n)$  \\

 $Mor= [(X,n),(Y,m)]= \left\{ \begin{array}{lll}
 \lim\limits_{\stackrel{N\to\infty}{N-n\,\in\mathbb{Z}}} \left[ \Sigma^{N-n}X,\Sigma^{N-m}Y \right]  & if & m-n\in\mathbb{Z} \\
 0 & if & m-n \notin\mathbb{Z}
 \end{array}\right. $\\
 
We define $\pin$-equivariant suspension spectra similarly. For our purposes, it suffices to use the  following real irreducible representations of $\pin$: \\
 
 $\left\{ \begin{array}{cl}
 \mathbb{R} & \text {with trivial action;}\\
 \widetilde{\mathbb{R}} & \left\{ \begin{array}{ll}
 j & \text{ acts multiplication by } -1,\\
 S^1 & \text{ acts trivially};
 \end{array} \right. \\
 \mathbb{H} & \text{action by } \pin \text{ via left multiplication.}
 \end{array} \right.$ \\
 
 Note that
 $$\pin \subset SU(2) = S(\mathbb{H}).$$
 
 We define a $\pin$-equivariant suspension spectrum as a quadruple $(X,n_{\mathbb{R}},n_{\widetilde{\mathbb{R}}},n_{\mathbb{H}})$.  

 In our case, the finite dimensional approximation to $V$ decomposes as
$$\begin{array}{ccccc}
V_\lambda^\mu & \cong & \widetilde{\mathbb{R}}^a & \oplus & \mathbb{H}^b \\
& & \uparrow &  & \uparrow \\
& & \text{forms} &  & \text{spinors}
 \end{array}$$
 
 We define the {\em Seiberg-Witten Floer (equivariant) spectrum} of $Y$ to be 
 $$\swf(Y):=\,\Sigma^{\mathbb{H}\frac{n(Y,g)}{4}} \,\Sigma^{-V_{\lambda}^{0}} \,I_{\lambda}^{\mu}.$$ 
 Then, we have
  $$H^{\pin}(\swf(Y); \mathbb{F}) = \swfh^{\pin}(Y; \mathbb{F}).$$
Other theories can be obtained by applying various generalized homology functors to $\swf(Y)$.

 \section{Homology cobordism invariants}
 \label{sec:abc}
 
We plan to use $\swfh^{\pin}(Y; \mathbb{F})$ to construct a map  $\beta:\Theta_3^H\to\mathbb{Z}$ satisfying the properties
\begin{itemize}
\item $\beta(-Y)=-\beta(Y)$;
\item $\beta(Y)(mod2)=\mu(Y).$
\end{itemize}
In the process, we will also obtain two other maps $\alpha, \gamma : \Theta_3^H\to\mathbb{Z}$.
 
\subsection{The module structure on equivariant homology} 
Recall that if we have a Lie group $G$ acting on a space $X$, the Borel homology
$$H^G_*(X)= H_*(X\times_G EG)$$ is a module over $H_G^*(pt)=H^*(BG)$. In our setting, we take $G = \pin$, and we are interested in understanding $H^*(B\pin; \mathbb{F})$. 

We have a fibration  
 $$\begin{tikzcd}
\pin\arrow{r}{i}  & SU(2) \arrow{d}{\psi}  \\
 {} &\mathbb{RP}^2
 \end{tikzcd}  $$
where $i$ is the inclusion and $\psi$ is the composition of the Hopf fibration with the involution on $S^2$. This fibration gives another fibration: 
  $$\begin{tikzcd}
\mathbb{RP}^2 \arrow{r}  & B\,\pin  \arrow{d}   \\
  {} & B\,SU(2)=\mathbb{HP}^\infty
  \end{tikzcd} $$

 The cohomology of $\mathbb{RP}^2$ in degrees  $0,\:1,\: 2$  is  
$$\includegraphics{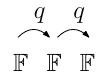}$$
\medskip

The cohomology of $B\,SU(2)$ in degrees  $0,\:1,\: 2,\:3,\:\dots$  is
$$\includegraphics{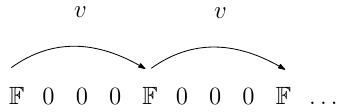}$$ 

\vspace{.2cm}

The Leray-Serre spectral sequence associated to the above fibration has no room for higher differentials. Thus, the cohomology groups of $B\pin$ in degrees $0,\:1,\: 2, \:3,\: 4, \:\dots$ are
$$\includegraphics{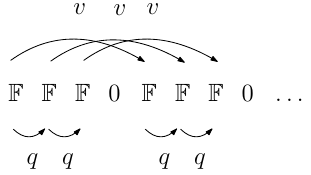}$$

The multiplicative property of spectral sequence gives a ring isomorphism
$$ H^*(B\,\pin;\mathbb{F})\cong \mathbb{F}[q,v]/(q^3), \: \text{deg}(v)=4, \: \text{deg}(q)=1. $$

Thus, if we have a space $X$ with a $\pin$-action, its Borel homology has an action by the ring above, with $q$ and $v$ decreasing grading by $1$ and $4$, respectively.

\subsection{Three infinite towers}
Recall that the Seiberg-Witten Floer spectrum $\swf(Y)$ is a formal (de-)suspension of the Conley index $I_\lambda^\mu$. This latter space is a finite $\pin$-equivariant $CW$-complex.

Let $(I_\lambda^\mu)^{S^1}$ denote the fixed points set of $I_\lambda^\mu$ under the action of the subgroup $S^1 \subset \pin$. Note that $(I_\lambda^\mu)^{S^1}$ picks up the part of the flow that lives in the reducible locus $\{(a,\phi)|\phi=0\}$. Note that the Seiberg-Witten equations
$$ \left\{\begin{array}{c l} 
*da &=\rho^{-1} ((\phi\otimes\phi^*)_{\circ}) \\
\partiali\phi+\rho(a)\phi&=0 \end{array}\right.$$
produce a linear flow (given by $*da$) on the locus where $\phi=0$. In view of this, one can check that 
$$(I_\lambda^\mu)^{S^1} = S^{dim \; V^0_{\lambda}}.$$

Recall that we defined 
$$\swf(Y)= \Sigma^{\mathbb{H}\frac{n(Y,g)}{4}} \: \Sigma^{-V_\lambda^0} \: I_\lambda^\mu.$$

Therefore,
$$(\swf(Y))^{S^1}=S^{n(Y,g)}.$$ 
Intuitively, $\swf(Y)$ is made of a reducible part $S^{n(Y,g)}$  and some free cells as the irreducible part. 
$$(\text{sphere}) \subset \swf(Y) \longrightarrow \swf(Y) /{\text{sphere}} \: \: \rotatebox[origin=c]{90}{$\circlearrowleft$} \: \: \pin \text{  acts freely}.$$
\smallskip

The $\pin$-equivariant Seiberg-Witten Floer homology
$$\swfh^{\pin}_*(Y; \mathbb{F})= \widetilde{H}_*^{\pin}(\swf(Y);\mathbb{F})$$ 
is a module over $\mathbb{F}[q,v]/(q^3)$, where $\mathbb{F}=\mathbb{Z}/2$.

There is a \textbf{localization theorem} in equivariant cohomology, which gives
$$V^{-1}\widetilde{H}^*_{\pin}(\swf(Y); \mathbb{F})= V^{-1}\widetilde{H}^*_{\pin}(S^{n(Y,g)}; \mathbb{F}).$$

Note that $ \widetilde{H}^*_{\pin}(S^{n(Y,g)}; \mathbb{F}) = H^{*\,-n(Y, g)}(B\pin; \mathbb{F})$. 

We can re-interpret the localization theorem in terms of Borel homology rather than Borel cohomology. Since we work over the field $\mathbb{F}$, Borel homology is simply the dual space to Borel cohomology, in each grading. Further, we can recover the module action on Borel homology from the one on Borel cohomology. The upshot is that $\swfh^{\pin}_*$ is the homology of a complex (the equivariant cellular complex of the Conley index) of the form:
\begin{center}
\includegraphics{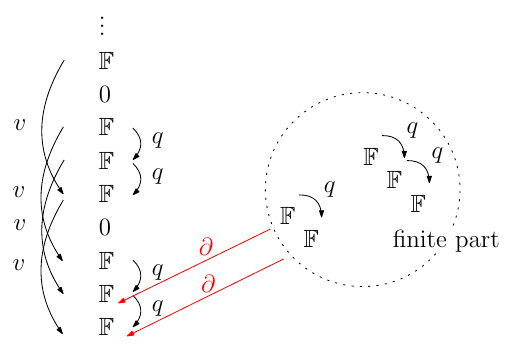}
\end{center}

The finite part can be any finite dimensional vector space, with some $v$ and $q$ actions, and with various differentials $\del$, either inside of it, or relating it to the infinite towers. In any case, there are always three infinite $v$ towers in the chain complex, and they produce three such towers in homology. The towers correspond to the $S^1$-fixed point set of $\swf(Y)$, and the finite part comes from the free cells.

Because $(\swf(Y))^{S^1}=S^{n(Y,g)}$ and $n(Y, g) \equiv 2\mu$ (mod $4$), we see that

\begin{itemize}
\item In the $1^{st}$ tower, all elements are in degree 
$ 2\mu\: (\text{mod }4);$

\item In the $2^{nd}$ tower, all elements are in degree
$ 2\mu +1 \: (\text{mod }4);$

\item In the $3^{rd}$ tower, all elements are in degree
$2\mu +2 \: (\text{mod }4).$
\end{itemize}

\subsection{Definition of the invariants}
 Let the lowest degrees of each infinite $v$-tower in homology be $A,\,B,\,C \in \mathbb{Z}$. 

Let $\alpha,\,\beta,\,\gamma\, \in \mathbb{Z}$ be invariants of $Y$ as follows:
$$\alpha= \frac{A}{2} ,\  \ \beta=\frac{B-1}{2} , \ \ \gamma=\frac{C-2}{2}.$$
Observe that
$$ \begin{array}{ccccccl}
 A & , & B-1 & , & C-2 & \equiv & 2\mu \ (\text{mod }4)\\
 & & & & & & \\
 \alpha= \frac{A}{2} & , & \beta=\frac{B-1}{2} & , & \gamma=\frac{C-2}{2} & \equiv &  \mu \ (\text{mod }2).
\end{array} $$

Furthermore, because of the module structure, we must have
$$\alpha \geq \beta \geq \gamma.
$$

\subsection{Descent to homology cobordism}
Next, we will check that $\alpha,\: \beta , \: \gamma$ descend to maps $\Theta_3^H \to \mathbb{Z}$. This uses the construction of cobordism maps on Seiberg-Witten Floer spectra from \cite{man03}.

Let $W^4$ be a smooth oriented $spin(4)$ cobordism with $b_1(W)=0$, with $\del W = (-Y_0) \cup Y_1$. (For our purposes, we are  interested in the case when $W$ is a homology cobordism  between homology $3$-spheres $Y_0$ and $Y_1$.) Basically, one can look at the $SW$-equations on $W$ and do a finite dimensional approximation to the solution space. This is somewhat similar to what we do in the $3$-dimensional  case. There is some more work to be done for cobordisms, but here we skip the details. The final result is a stable equivariant map between two suspension spectra:
$$\Psi_W : \Sigma^{m\mathbb{H}}\swf(Y_0) \to \Sigma^{n\widetilde{\mathbb{R}}}\swf(Y_1),$$

Here $m\mathbb{H}$ is the direct sum of the $m$ copies of quaternionic representation and similarly $n\widetilde{\mathbb{R}}$ is the direct sum of the $n$ copies of the sign representation.  Also, 
$$m=\frac{-\sigma(W)}{8} = ind\Di, \ \ \ n=b_2^+(W)=ind(d^+).$$

\begin{example} 
Assume $Y_0=Y_1=S^3$. Then, we can fill in the cobordism $W$ with two copies of $B^4$ and get a closed four-manifold $X$. We have $\swf(S^3)=S^0$, so in this case we get  a map 
$$\Psi_W: (m\mathbb{H})^+\to (n\widetilde{\mathbb{R}})^+ $$
This map is the Bauer-Furuta invariant of $X$, a stable homotopy refinement of the Seiberg-Witten invariant. See \cite{bf04}.
\end{example}

Now, suppose $W$ is a smooth oriented homology cobordism between homology spheres $Y_0$ and $Y_1$. There is a unique $spin(4)$ structure on $W$. Moreover, we have $b_1(W)=0$ and there is  $m=n=0$. Let $F_W$ be the homomorphism induced on $\pin$-equivariant homology by the map $\Psi_W$:

\begin{center}
\includegraphics{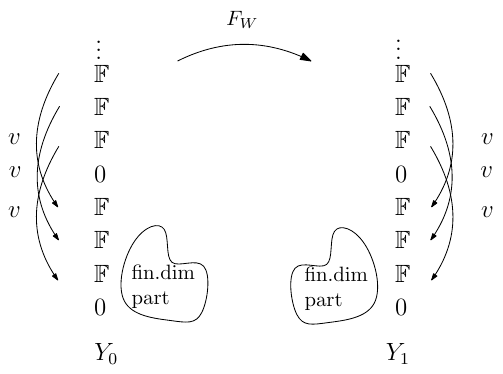}
\end{center}

It follows from equivariant localization that in  degrees $k\gg 0 $, the map $F_W$ is an isomorphism. Further, $F_W$ is a module map, so we have a commutative diagram
$$\begin{tikzcd}
\mathbb{F}\arrow{r}{F_W} \arrow{d}{v}  & \mathbb{F} \arrow{d}{v}  \\
 \mathbb{F} \arrow{r}{F_W} & \mathbb{F}
 \end{tikzcd}  $$
 
Because of the module structure, we cannot have $\alpha(Y_1)<\alpha(Y_0)$, and likewise for $\beta$ and $\gamma$. In conclusion,
 $$\begin{array}{ccl}
 \alpha(Y_1)  & \geq & \alpha(Y_0),\\
 \beta(Y_1)  & \geq & \beta(Y_0),\\
  \gamma(Y_1)  & \geq & \gamma(Y_0).
 \end{array}$$
 
On the other hand, if we reverse the orientation and the direction of $W$ we get a homology cobordism from $Y_1$ to $Y_0$. By the same logic, we get
  $$\begin{array}{ccl}
  \alpha(Y_0)  & \geq & \alpha(Y_1),\\
  \beta(Y_0)  & \geq & \beta(Y_1),\\
   \gamma(Y_0)  & \geq & \gamma(Y_1).
  \end{array}$$
  
Thus, we have equalities, and obtain the following corollary.

\begin{corollary}
The invariants $\alpha,\: \beta , \: \gamma$ descend to maps $\Theta_3 ^H \to \mathbb{Z}$.
\end{corollary}

\section{Duality}
\label{sec:duality}

Recall from Section~\ref{sec:strategy} that, in order to prove Theorem~\ref{thm:man} and hence Theorem~\ref{thm:tc}, it suffices to construct an invariant $\beta$ with the three properties listed there. We have already checked that our $\beta$ reduces mod $2$ to the Rokhlin invariant, and also that it descends to $\Theta_3^H$. It remains to prove that
$$ \beta(-Y) = - \beta(Y).$$
We will do this in this section. In the process, we will also find that the other two invariants, $\alpha$ and $\gamma$, satisfy  
$$\alpha(-Y)=-\gamma(Y), \ \ \gamma(-Y)=-\alpha(Y).$$

\subsection{Orientation reversal}
Consider a homology $3$-sphere $(Y,g)$ and change its orientation to get $(-Y,g)$. This changes the direction of the $SW$ flow equation.
$$\dot{x}=-SW(x(t)) \hspace{.5cm} \rightsquigarrow \hspace{.5cm} \dot{x}=SW(x(t)) $$

In the finite dimensional approximation $V_\lambda^\mu$, we can choose index pairs for the forward and reverse flows, $(N,L_{+})$ and $(N,L_{-})$, such that $N$ is a codimension $0$ submanifold (with boundary) of $V_\lambda^\mu$, and
$$L_{+} \cup L_{-}= \partial N \:,\: \partial L_{+} = \partial L_{-} = L_{+} \cap L_{-}.$$

Non-equivariantly we have a duality
\begin{equation}
\label{eq:ndual}
\widetilde{H}_*(N/L_{+}) = \widetilde{H}^{d-*}(N/L_{-}),
\end{equation}
where $d=dim \; V_\lambda^\mu = dim \; N$. This follows from Alexander duality. Indeed, one can find an embedded $X\subset V_\lambda^\mu \times \mathbb{R} = \mathbb{R}^{d+1}$ such that   $$X\simeq N/L_{+} \hspace{.5cm} \text{and} \hspace{.5cm} \mathbb{R}^{d+1} - X\simeq N/L_{-} $$ 

If you embed a space in Euclidian space then the homology of the space and the cohomology of the complement are related by Alexander duality. Equation \eqref{eq:ndual} follows from here.

\subsection{Spanier-Whitehead duality} We seek an equivariant analogue of \eqref{eq:ndual}. Before getting to that, it is helpful to first understand the stable homotopy version of \eqref{eq:ndual},  which is called Spanier-Whitehead duality.

Non-equivariantly, consider a suspension spectrum, the  formal (de-)suspension of a topological space $X$:
 $$\mathcal{Z} = (X, k) = \Sigma^{-k}X.$$

Embed $X \hookrightarrow S^N$ for some $N\gg 0$.  

\begin{definition} The \emph{Spanier-Whitehead dual} of $\Sigma^{-k}X $ is 
$$D(\Sigma^{-k}X) := \Sigma^{k}(\Sigma^{-(N-1)}(S^N-X)).$$
\end{definition}

\begin{example}
Consider  $V=\mathbb{R}^2$ and embed $S^1 \hookrightarrow S^2$. Then, the complement is homotopy equivalent to $S^0$. In this case, $N=2$, and $D(S^1)= \Sigma^{-1}S^0 = S^{-1}$.
\smallskip
\begin{center}
\includegraphics{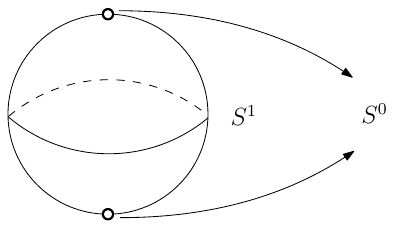}
\end{center}
\end{example}
\smallskip

More generally, we have $D(S^k)= S^{-k}= (S^0,k)$. Furthermore, the dual of a wedge product of spaces is the wedge product of the duals, and similarly for smash products.

It follows from Alexander duality that, if $\mathcal{Z}$ is a suspension spectrum,
$$\widetilde{H}_k(\mathcal{Z})=\widetilde{H}^{-k}(D(\mathcal{Z})).$$ 

Similar equalities also hold for other generalized homology theories. \\

Here is the equivariant analogue:

\begin{definition}
Let $G$ be a Lie group, $X$ a $G$-space, and $W$ a representation of $G$. Let us embed $X \hookrightarrow V^+$, for some representation $V$ of $G$. We define the {\em equivariant Spanier-Whitehead dual} of a formal de-suspension $\Sigma^{-W}X$ by:
$$ D(\Sigma^{-W}X) := \Sigma^{W}(\Sigma^{-V} \Sigma^\mathbb{R}(V^+-X)). $$ 
\end{definition}

Given our choice of the index pairs for the forward and reverse flow in $V^{\mu}_{\lambda}$, we find that the Seiberg-Witten Floer spectra of $Y$ and $-Y$ are related by $\pin$-equivariant duality:
 $$D(\swf(Y))=\swf(-Y).$$

For example, consider the Poincare homology sphere $Y=\Sigma(2,3,5)$. We have
$$\begin{array}{ccccc}
\swf(Y) & = & S^2 & = & \Sigma^{\frac{1}{2}\mathbb{H}} S^0, \\
 & & & & \\
\swf(-Y) & = & S^{-2} & = & \Sigma^{\frac{-1}{2}\mathbb{H}} S^0.
\end{array}$$

\subsection{Duality for equivariant homologies}  
Given a $G$-equivariant suspension spectrum $\mathcal{Z}$, we seek to investigate the relation between $$\widetilde{H}_*^G(\mathcal{Z})\ \ \text{and} \ \ \widetilde{H}_G^{-*}(D(\mathcal{Z})).$$ In the non-equivariant case, they are the same. However, this cannot be true equivariantly! Indeed, usually the left hand side is infinite in the positive direction only (that is, non-trivial in some degrees $k\gg 0$, and trivial for $k \ll 0$), and the right hand side is infinite in the negative direction only (that is, nontrivial in some degrees $k\ll 0$, but trivial for $k \gg 0$). \\

Recall the definitions of the Borel cohomology and homology, respectively: 
$$\begin{array}{ccc}
H_G^*(X) & = & H^*(X\times_G EG), \\
& & \\
H^G_*(X) & = & H_*(X\times_G EG).\\
\end{array}$$

\begin{definition}
The {\em co-Borel homology} of an equivariant suspension spectrum is defined as 
$$c\widetilde{H}_*^G(\mathcal{Z}) = \widetilde{H}_G^{-*}(D(\mathcal{Z})),$$ 
where $\mathcal{Z}= \Sigma^{-V}X$ for some $X$. 
\end{definition}
We want to understand relation between Borel and co-Borel homology. This relation goes through {\em Tate homology}.

\begin{definition}
Following \cite{gm95}, we let the {\em Tate homology} of $\mathcal{Z}= \Sigma^{-V}X$ be 
$$t\widetilde{H}_*^G(\mathcal{Z})= c\widetilde{H}_*^{G} (\widetilde{EG}\wedge \mathcal{Z}),$$
 where $\widetilde{EG}$ is the unreduced suspension of the $EG$.
\end{definition}

The main property of Tate homology that we need is that 
$$t\widetilde{H}_*^G(\mathcal{Z})= 0 \text{ if } G \text{ acts freely on } \mathcal{Z}.$$

Also, the Borel, co-Borel and Tate homologies  ($\widetilde{H}_*^G$, $c\widetilde{H}_*^G$, and $t\widetilde{H}_*^G$) satisfy the usual excision and suspension properties of homology.

The relation between these homologies that we alluded to is given by the {\em Tate-Swan exact sequence}:
$$\dots\to \widetilde{H}^G_{n-dimG}(\mathcal{Z}) \to c\widetilde{H}^G_{n}(\mathcal{Z}) \to t\widetilde{H}^G_{n}(\mathcal{Z}) \to \widetilde{H}^G_{n-dimG-1}(\mathcal{Z}) \to \dots $$

\begin{example}
The simplest example is when $G= 1$. A point acts freely on any space, so $t\widetilde{H}^G_* = 0 $, and the Tate-Swan exact sequence gives rise to the usual Alexander duality:
$$\widetilde{H}^G_*(\mathcal{Z}) = c\widetilde{H}^G_*(\mathcal{Z}) = \widetilde{H}_G^{-*}(D(\mathcal{Z})).$$
\end{example}

\begin{example}
Let $G=S^1$ and $\mathcal{Z}=\swf(Y)$ for a homology sphere $Y$. Then,
$$t\widetilde{H}^{S^1}_*(\mathcal{Z}) = t\widetilde{H}^{S^1}_*(\text{fixed point set}) = t\widetilde{H}_*^{S^1}(\text{sphere}) = \mathbb{Z}[U, U^{-1}].$$
The Tate-Swan exact sequence looks like:
\begin{center}
\includegraphics{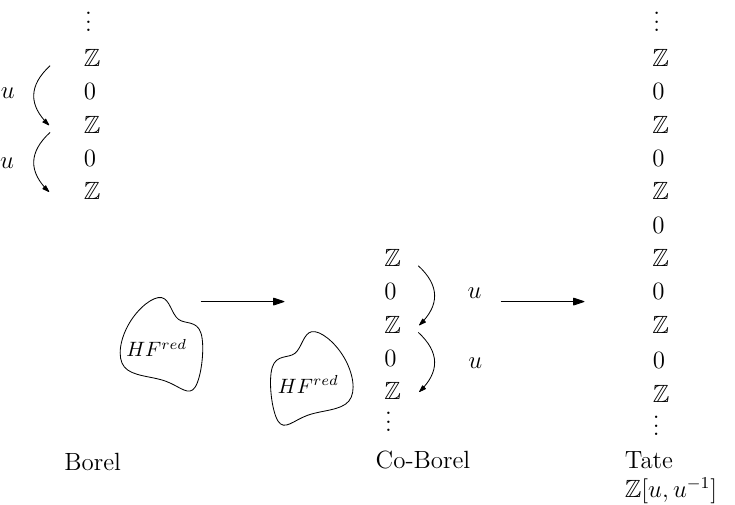}
\end{center}
\end{example}
\smallskip

Here is a dictionary between the different $S^1$-equivariant theories, and their counterparts in Heegaard Floer homology, as well as in the Kronheimer-Mrokwa version of monopole Floer homology:
$$\begin{array}{ccc}
\begin{array}{lcc}
 &  & \\
  &  & \\
  \text{ Borel } & \widetilde{H}_*^{S^1} & \swf(Y) \\
  \text{ co-Borel } & c\widetilde{H}_*^{S^1} & \swf(Y)\\
     \text{ Tate } & t\widetilde{H}_*^{S^1} & \swf(Y)\\
     \text{ Non-equivariant } & \widetilde{H}_* & \swf(Y) \\
      &  & 
\end{array} & \left| \begin{array}{clc}
    & \text{ Heegaard Floer} & \\
     &  & \\
 & HF^+ & \\
  & HF^- & \\
  & HF^\infty & \\
  &\widehat{HF} & \\
       &  & 
\end{array} \right| & \begin{array}{clc}
 
 & \text{Kronheimer-Mrowka} & \\
  &  & \\
 & \widecheck{HM} & \\
 & \widehat{HM} & \\
 & \overline{HM} & \\
 & \widetilde{HM} & \\
      &  & 
\end{array}

\end{array}$$

Now we arrive at the case of interest to us. Let $G=\pin$ and $\mathcal{Z}=\swf(Y)$. The Tate homology is
$$t\widetilde{H}_*^{\pin} (\swf(Y); \mathbb{F}) =  t\widetilde{H}^{S^1}_*(S^1\text{-fixed point set}; \mathbb{F}) = \mathbb{F}[q, v, v^{-1}]/(q^3).$$

This is related to the Borel homology $\widetilde{H}_*^{\pin} (\swf(Y); \mathbb{F})$ and the co-Borel homology $c\widetilde{H}_*^{\pin} (\swf(Y); \mathbb{F}) = \widetilde{H}^{-*}_{\pin} (\swf(-Y); \mathbb{F})$ by a Tate-Swan sequence of the form
\begin{center}
\includegraphics{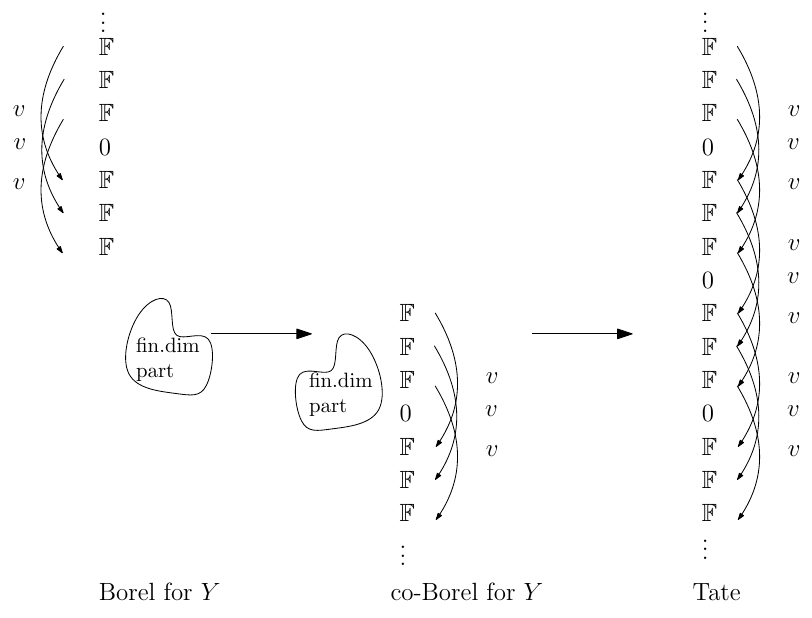}
\end{center}
\bigskip

\subsection{Conclusion}
By analyzing the Tate-Swan exact sequence above, we see that the position of the bottom elements  in the three infinite towers in Borel homology determines the position of the top elements in the three towers in co-Borel homology, and hence that of the bottom elements in the three towers  in the Borel homology of $D(\swf(Y)) = \swf(-Y)$. In this process,  we take the negative of the grading  when we pass from co-Borel homology to Borel cohomology of the dual. Hence, the first tower in $\widetilde{H}_*^{\pin} (\swf(Y); \mathbb{F})$ corresponds to the third tower in $\widetilde{H}_*^{\pin} (\swf(-Y); \mathbb{F})$, the second tower corresponds to the second, and the third to the first. Therefore, we obtain 
$$\begin{array}{ccc}
\gamma(-Y) & = & -\alpha(Y), \\
\beta(-Y) & = & -\beta(Y), \\
\alpha(-Y) & = & -\gamma(Y).
\end{array}$$

We have now established all the desired properties of $\beta$:
$$\left. \begin{array}{c}
\beta: \Theta_3 ^H \to \mathbb{Z}\\
\begin{array}{ccc}
\beta(Y)  & \equiv & \mu(Y) \ (\text{mod }2 )  \\
\beta(-Y) & = & -\beta(Y)
\end{array}    
\end{array} \right\rbrace \Longrightarrow $$
(by the discussion in Section~\ref{sec:strategy}) the short exact sequence 
$$0\longrightarrow Ker\mu\longrightarrow \Theta_3 ^H \stackrel{\mu}{\longrightarrow} \mathbb{Z}/2\longrightarrow 0$$
does not split. 

Combining this with the work of Galewski-Stern and Matumoto (see Section~\ref{sec:manifolds}), we complete the proof of Theorem~\ref{thm:tc}: For all $n\geq 5$, there exist non-triangulable manifolds of dimension $n$.

\section{Involutive Heegaard Floer Homology} 
\label{sec:hfi}
In this section we describe some recent joint work of Hendricks and the author \cite{hm15}.

\subsection{Motivation and outline}
We have constructed $\alpha, \: \beta, \: \gamma : \: \Theta_3^H \to \mathbb{Z}$, using $\pin$-equivariant Seiberg-Witten Floer homology. The issue is that Seiberg-Witten Floer homology is rather hard to compute. There are some calculations:
\begin{itemize}
\item  For Seifert fibrations by M.Stoffregen  \cite{sto15};
\item For surgeries on alternating knots by F.Lin \cite{lin15}.
\end{itemize}

However, most of these computations are based on the isomorphism between monopole Floer homology and Heegaard Floer homology $\widecheck{HM} \cong {HF}^+$, since the latter is the more computable theory. 

Let us recall from \cite{os04a, os04b} that, to a homology $3$-sphere $Y^3$ with $Spin^c$ structure $\mathfrak{s}\in Spin^c(Y)$, Ozsv\'ath and Szab\'o assign the (plus version of) Heegaard Floer homology
$$HF^+(Y,\mathfrak{s}).$$
This is a module over $\mathbb{F}[U]$, and it corresponds to $S^1$-equivariant Seiberg-Witten Floer homology. Furthermore, there is a conjugation symmetry $\iota_*$ on Heegaard Floer homology, that takes the $spin^c$ structure to the conjugate $spin^c$ structure $\mathfrak{s} = \bar{\mathfrak{s}}$:
$$ HF^+(Y,\mathfrak{s}) \cong HF^+(Y,\bar{\mathfrak{s}}).$$

Recall also that $$\pin=S^1\cup jS^1.$$

Ideally, we would want to use $ \iota_* $ to construct a $\mathbb{Z}/2$-equivariant $HF^+$, which  should correspond to $\pin$-equivariant Seiberg-Witten Floer homology. However, as we shall explain later, this is beyond current Heegaard Floer technology. 

Instead, we define a theory denoted
$$HFI^+(Y,\mathfrak{s}),$$  and called {\em involutive Heegaard Floer homology}. This  corresponds to $\mathbb{Z}/4$-equivariant Seiberg-Witten Floer homology, where $\mathbb{Z}/4$ is the subgroup 
$$\mathbb{Z}/4 = \langle j \rangle \subset \pin.$$
The involutive theory can be defined for all (pairs of) $Spin^c$ structures, but it only contains new information (compared to $HF^+$) for self-conjugate $Spin^c$ structures, i.e., those that come from spin structures.

\begin{theorem} [\cite{hm15}] Let $Y$ be a $3$-manifold and $\mathfrak{s}\in Spin(Y)$. Then, the isomorphism class of the $HFI^+(Y,\mathfrak{s})$, as a module over the cohomology ring $H^*(B\mathbb{Z}/4; \mathbb{F}) = \mathbb{F}[Q,U]/(Q^2)$, $(deg(U)=-2, \: deg(Q)=-1 )$ is an invariant of $(Y,\mathfrak{s})$. 
\end{theorem} 

The invariant $HFI^+$ does not have the full power of $\swfh^{\pin}$, but it is more computable. Indeed, in principle $HFI^+$ is algorithmically computable for large surgeries on all knots $K\subset S^3$, by using grid diagrams. More explicitly, one can get concrete formulas for large surgeries on $L$-space knots and quasi-alternating knots. (Here, ``large'' means with surgery coefficient an integer greater or equal than the genus of the knot.)

\subsection{Homology cobordism invariants} 
From the different Floer homologies, one can construct various homology cobordism invariants. For $Y\in \mathbb{Q}HS^3$ and $\mathfrak{s}\in Spin^c(Y)$,  the minimal grading of the infinite $U$-tower in $HF^+(Y,\mathfrak{s})$ is called the {\em Ozsv\'ath-Szab\'o correction term} 
$$d(Y, \mathfrak{s}) \in \mathbb{Q}.$$ This was defined in \cite{os03}, and is the analogue of (twice) the Fr{\o}yshov invariant $\delta$ from monopole Floer homology \cite{fro10}. When $Y$ is a $\mathbb{Z}/2$-homology sphere, it has a unique $Spin^c$ structure $\mathfrak{s}$. Further, if $Y$ is a $\mathbb{Z}$-homology sphere, then $d$ takes even integer values. Thus, we get  homomorphisms 
$$ \begin{array}{cccc}
d: & \Theta^3_{\mathbb{Z}/2} & \rightarrow & \mathbb{Q}\\
   & & & \\
   & \big \uparrow &  & \bigcup \\
   & & & \\
   & \Theta^3_{\mathbb{Z}} & \rightarrow & 2\mathbb{Z}
\end{array}$$
Here, we changed notation and let $\Theta^3_{\mathbb{Z}}$ be the homology cobordism group with integer coefficients, previously denoted by $\Theta_3^H$. This allows us to distinguish it from $\Theta^3_{\mathbb{Z}/2}$, the homology cobordism group with $\mathbb{Z}/2$ coefficients. The latter is generated by $\mathbb{Z}/2$-homology spheres, and uses a weaker equivalence relation, given by the existence of a cobordism $W$ with $H_*(W, Y_i; \mathbb{Z}/2)=0$.

From $HFI^+$, if $Y\in \mathbb{Q}HS^3$ and $\mathfrak{s}=\widebar{\mathfrak{s}}$ then, in a similar manner, we get two new invariants
$$ \underline{d}(Y,\mathfrak{s}) , \: \widebar{d}(Y,\mathfrak{s}) \in \mathbb{Q}. $$

These descend to maps
$$ \begin{array}{ccccc}
\underline{d}, \: \widebar{d}: & \Theta^3_{\mathbb{Z}/2} &  \rightarrow & \mathbb{Q}  \\
   & & &  \\
   & \big \uparrow &  & \bigcup  \\
   & & & \\
   & \Theta^3_{\mathbb{Z}} & \rightarrow & \mathbb{Z} 
\end{array}$$
One should mention here that $\widebar{d}$ and $\underline{d}$ are not homomorphisms.

Here is a list of the different homology cobordism invariants that we mentioned:
$$\begin{array}{clcl}
\text{from } & \swfh^{S^1} \ (\text{or } HF^+) & \stackrel{get}{\rightsquigarrow}& \delta \ \ (\text{or } d); \\
\text{from } & \swfh^{\pin}& \stackrel{get}{\rightsquigarrow}& \alpha, \beta, \gamma ; \\
\text{from } & \swfh^{\mathbb{Z}/4} \ (\text{or } HFI^+) & \stackrel{get}{\rightsquigarrow}& \underline{d} , \widebar{d} .
\end{array}$$

\begin{example} For the Brieskorn sphere $\Sigma(2,3,7)$, we have 
$$\widebar{d}(\Sigma(2,3,7)) = d((\Sigma(2,3,7))) = 0\ \text{ but } \ \underline{d}(\Sigma(2,3,7))= -2.$$
\end{example} 

On the other hand, $d(Y) = \underline{d}(Y) = \widebar{d}(Y)$ if $Y$ is an $L$-space, e.g., a double branched cover of an alternating knot.

Thus, we obtain the following application of involutive Heegaard Floer homology:
\begin{corollary}
$\Sigma(2,3,7)$ is not $\mathbb{Z}/2$-homology cobordant to any $L$-space.
\end{corollary}

The same result can be obtained using $\alpha$, $\beta$, $\gamma$ from $\swfh^{\pin}$. See Corollary~\ref{cor:app} below for another application of $HFI^+$, for which we do not yet have a proof based on Seiberg-Witten theory.

Let us end by noting that $\underline{d}, \: \widebar{d} \equiv d \  (\text{mod }2)$, but (unlike in the case of $\alpha$ , $\beta$, $\gamma$), in general $$\underline{d}, \: \widebar{d} \nequiv \mu \  (\text{mod }2).$$ Thus, we cannot re-disprove the triangulation conjecture by using  $\underline{d}, \: \widebar{d}$. On the other hand, $\underline{d}, \: \widebar{d}$ are more computable than $\alpha$, $\beta$, $\gamma$.

\subsection{Concordance Invariants}
Whenever we have an invariant of homology cobordism, we get invariants of smooth knot concordance by either doing surgeries on the knot, or by taking the double branched cover. Let us focus on surgeries.

In Heegaard Floer homology, there is an invariant of knots $K\subset S^3$:
$$V_0(K) =  \frac{p-1}{8} -\frac{1}{2} d(S^3_{p}(K)), \ p \in \mathbb{Z}, \ p > 0.$$

If $K$ is smoothly concordant to $K'$ ($K\sim K'$), then $S^3_p(K)$ is homology cobordant to $S^3_p(K')$, and hence $V_0(K) = V_0(K')$.

Similarly, in the involutive theory we can define
$$\Vl_0(K) = \frac{p-1}{8}  -\frac{1}{2} \underline{d}(S^3_{p}(K)),$$
$$\overline{V}_0(K) = \frac{p-1}{8}  -\frac{1}{2} \widebar{d}(S^3_{p}(K))$$
for integers $p \geq g(K)$, where $g(K)$ is the genus of the knot. (Conjecturally, the same formulas hold for all integers $p > 0$.)

For instance, $\Sigma(2,3,7)$ is $+1$-surgery on the figure-eight knot $4_1$ of genus one:  
\begin{center}
\includegraphics{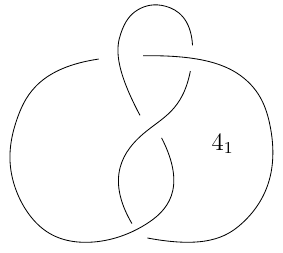}
\end{center}

We get
$$ V_0 (4_1) = \overline{V}_0(4_1) = 0, \:  \Vl_0(4_1) = 1. $$ 

Thus, $\Vl_0$ detects the non-sliceness of $4_1$, unlike $V_0$ (or other similar invariants from Khovanov and Floer theories). Alternate proofs that $4_1$ is not slice can be given using the Fox-Milnor condition on the Alexander polynomial $\Delta_K$, or using the Rokhlin invariant of surgeries. 

The knot invariants $V_0, \overline{V}_0$ and $ \Vl_0$ can be calculated explicitly for $L$-space and quasi-alternating knots. From here, we get constrains on which $3$-manifolds can be homology cobordant to other $3$-manifolds. Here is a sample application:
\begin{corollary}
\label{cor:app}
Let $K$ and $K'$ be alternating knots such that $\sigma(K) \equiv 4\cdot \operatorname{Arf}(K) + 4 \pmod 8$. If $S^3_{p}(K)$ and $S^3_{p}(K')$ are $\Z_2$-homology cobordant for some odd $p \geq \max(g(K), g(K'))$, then $\sigma(K) =\sigma(K')$. (Here, $\sigma$, $\operatorname{Arf}$, and $g$ denote the signature, the Arf invariant, and the Seifert genus of a knot, respectively.) 
\end{corollary}

\subsection{Construction}
We now sketch the construction of $HFI^+$.

First, recall from \cite{os04a} that the Heegaard Floer homology $HF^+$ is computed from a Heegaard diagram associated to a closed, connected, oriented three-manifold $Y$:
$$\mathcal{H}= (\Sigma, \overrightarrow{\alpha} , \overrightarrow{\beta}, z )$$ 
Here, $\overrightarrow{\alpha} , \overrightarrow{\beta}$ are some collections of curves that determine the two handlebodies in the Heegaard splitting. 

Out of this, by doing Lagrangian Floer homology on the symmetric product of the Heegaard surface $\Sigma$, Ozsv\'ath and Szab\'o get a chain complex 
$$CF^+(\mathcal{H}, \mathfrak{s}).$$

Ozsv\'ath and Szab\'o prove that for different Heegaard diagrams $\mathcal{H}, \mathcal{H'}$ of $Y$, there is a chain homotopy equivalence 
$$\Phi(\mathcal{H}, \mathcal{H'}) : CF^+(\mathcal{H}, \mathfrak{s}) \to CF^+(\mathcal{H'}, \mathfrak{s})$$

Therefore, the isomorphism class of $HF^+(\mathcal{H},\mathfrak{s})$, as a $\mathbb{F}[U]$-module, is an invariant of $(Y, \mathfrak{s})$; see \cite{os04a}.

There is a stronger statement proved by Juh\'asz and D.Thurston \cite{jk12}. They showed the naturality of the invariant. Naturality says that we can choose
$$\begin{array}{ccl}
\Phi(\mathcal{H}, \mathcal{H'}) & \text{ such that } & \left\lbrace \begin{array}{l}
\Phi(\mathcal{H}, \mathcal{H}) = Id \\
\Phi(\mathcal{H'}, \mathcal{H''}) \circ \Phi(\mathcal{H}, \mathcal{H'}) \sim  \Phi(\mathcal{H}, \mathcal{H'}) 
\end{array}   \right. 
\end{array}$$
where ``$\sim$'' denotes chain homotopy. This implies that the $\mathbb{F}[U]$-module
$$HF^+(\mathcal{H}, \mathfrak{s}) = HF^+(Y,\mathfrak{s}) $$ is a true invariant of $(Y,\mathfrak{s})$.

Next, let us consider the conjugation symmetry between Heegaard Floer complexes. There is an identification between the Heegaard diagrams
$$\mathcal{H}= (\Sigma, \overrightarrow{\alpha} , \overrightarrow{\beta}, z ) \stackrel{\cong}{\rightarrow} \bar{\mathcal{H} } = ( - \Sigma, \overrightarrow{\beta}, \overrightarrow{\alpha}, z ).$$

If $\mathfrak{s} = \bar{\mathfrak{s}}$ then $\iota: \: \:  CF^+(\mathcal{H},\mathfrak{s}) \to CF^+(\mathcal{H},\mathfrak{s}) $  is defined as composition of $\eta$ and $\Phi(\mathcal{H},\bar{\mathcal{H}})$:
$$
     \begin{tikzcd}[row sep=huge, column sep=huge, text height=1.5ex, text depth=0.25ex]
          \displaystyle CF^+(\mathcal{H},\mathfrak{s}) \arrow{r}{\eta}[swap]{\cong} \arrow[bend right=30]{rr}{\iota} & CF^+(\bar{\mathcal{H}},\bar{\mathfrak{s}}) \arrow{r}{\Phi(\mathcal{H},\bar{\mathcal{H}})} & CF^+(\mathcal{H},\bar{\mathfrak{s}})
     \end{tikzcd}$$

\begin{definition}
Whenever $\mathfrak{s} = \bar{\mathfrak{s}}$, we define 
$$ \begin{array}{rlcl}
 CFI^+(Y,\mathfrak{s})  = & \text{mapping cone of } \lgroup CF^+(\mathcal{H},\mathfrak{s}) &  \stackrel{1+\iota}{\longrightarrow} & CF^+(\mathcal{H},\mathfrak{s}) \rgroup \\
  = &  \text{mapping cone of } \lgroup  CF^+(\mathcal{H},\mathfrak{s}) & \stackrel{Q(1+\iota)}{\longrightarrow} & Q \cdot CF^+(\mathcal{H},\mathfrak{s}) [-1]  \rgroup
  \end{array}$$

\end{definition}

In the second row, we added a formal variable $Q$ of degree $-1$, and let $[-1]$ represent a shift in grading that cancels the shift by $Q$. This way, $CFI^+(Y,\mathfrak{s})$ becomes a module over the ring $\mathbb{F}[Q,U]/(Q^2)$. 

Using the naturality of $CF^+$, one can show that  $\Phi(\mathcal{H},\bar{\mathcal{H}} ) $ is well defined up to chain homotopy. From here it follows that the isomorphism class of $HFI^+(Y,\mathfrak{s})$ is a $3$-manifold invariant. 

\begin{remark}
We cannot show the naturality of $HFI^+(Y,\mathfrak{s})$ with the current methods. Naturality of $HF^+$ is used in the definition of $HFI^+$, and if we wanted naturality of $HFI^+$ we would need a kind of second order naturality for $HF^+$, which is not known. 
\end{remark}

\subsection{Intuition} Let us explain the motivation behind the definition of $CFI^+$ as the mapping cone of $1+\iota$.

If we wanted a Heegaard Floer analogue of $\pin$-equivariant Seiberg-Witten Floer homology, this would be a $\mathbb{Z}/2$-equivariant version of $HF^+$, with respect to the involution $\iota_*$.

A construction of $\mathbb{Z}/2$-equivariant Lagrangian Floer homology was proposed by Seidel and Smith in \cite{ss10}. It was inspired by $\mathbb{Z}/2$-equivariant  Morse theory. If we have a space $X$ wiht $\mathbb{Z}/2$ action then
$$H_*^{\mathbb{Z}/2} (X) = H (X\times_{\mathbb{Z}/2} E\mathbb{Z}/2 ),$$  
where $ X\times_{\mathbb{Z}/2} E\mathbb{Z}/2 $ is a bundle over $\mathbb{RP}^\infty = B\mathbb{Z}/2 $ with fibers $X$. Basically, we want to do Morse theory on this space. For this, take the standard Morse function on $\mathbb{RP}^N$ for large $N\gg 0$, and combine it with a family of Morse functions on the fibers, such that the critical points of the combined function lie over the critical points on $\mathbb{RP}^N$:

\begin{center}
\includegraphics{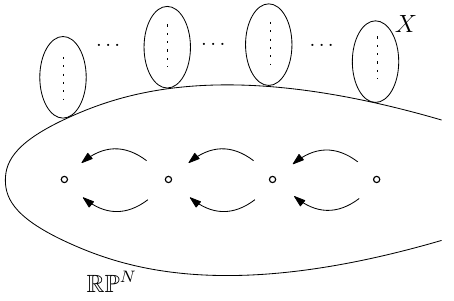}
\end{center}

We get a $\mathbb{Z}/2$-equivariant Morse complex of the form
$$C^{\mathbb{Z}/2}(X) =\xymatrixcolsep{3.5pc}
\xymatrix{
C(X)  &  C(X) \ar[l]_{Q\cdot (1+ \iota)} & C(X) \ar[l]_{Q\cdot (1+ \iota)}  \ar@/^1pc/[ll]^{Q^2 \cdot W}   & \dots \ar[l]_{Q\cdot (1+ \iota)}  \ar@/^1pc/[ll]^{Q^2 \cdot W}  \ar@/^3pc/[lll]^{Q^3 \cdot Z} 
}, $$
where $W$ is a chain homotopy between $\iota^2$ and the identity, $Z$ is a higher chain homotopy, and so on.

Thus, a $\pin$ version of Heegaard Floer homology should come from a chain complex of the form
$$CF^{\pin} = \xymatrixcolsep{3.5pc}
\xymatrix{
CF^+  &  CF^+ \ar[l]_{Q\cdot (1+ \iota)} & CF^+  \ar[l]_{Q\cdot (1+ \iota)}  \ar@/^1pc/[ll]^{Q^2 \cdot W}   & \dots \ar[l]_{Q\cdot (1+ \iota)}  \ar@/^1pc/[ll]^{Q^2 \cdot W}  \ar@/^3pc/[lll]^{Q^3 \cdot Z} 
}.$$

The problem is that, whereas to define $\iota$ we used naturality for $HF^+$, to define the higher homotopies in the infinite chain complex above we would need naturality to infinite order in Heegaard Floer theory. This seems intractable with current technology. Instead, we do a truncation, by setting 
 $$Q^2=0.$$

The truncated complex is the mapping cone $CFI^+$. This justifies the definition of $HFI^+$. To see that $HFI^+$ should correspond to $\mathbb{Z}/4$-equivariant Seiberg-Witten Floer homology, one needs the following algebraic-topologic fact.

\begin{lemma}
If we have a space $X$ with a $\pin$-action, then 
$$ H_*(C_*^{\pin} (X) /{Q^2}; \mathbb{F}) = H^{\mathbb{Z}/4}_*(X; \mathbb{F}).$$
\end{lemma}

We refer to \cite[Section 3]{hm15} for a proof.

\subsection{The structure of involutive Heegaard Floer homology} If $Y$ is a rational homology $3$-sphere equipped with a $Spin$ structure $\mathfrak{s}$, one can show that the involutive Heegaard Floer homology takes the following form: 
 \vskip.4cm
 \begin{center}
 \begin{picture}(0,0)%
\includegraphics{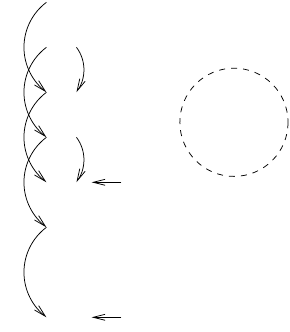}%
\end{picture}%
\setlength{\unitlength}{3158sp}%
\begingroup\makeatletter\ifx\SetFigFont\undefined%
\gdef\SetFigFont#1#2#3#4#5{%
  \reset@font\fontsize{#1}{#2pt}%
  \fontfamily{#3}\fontseries{#4}\fontshape{#5}%
  \selectfont}%
\fi\endgroup%
\begin{picture}(2889,3297)(1261,-3050)
\put(1911,157){\rotatebox{90.0}{\makebox(0,0)[lb]{\smash{{\SetFigFont{10}{12.0}{\rmdefault}{\mddefault}{\updefault}{\color[rgb]{0,0,0}$\dots$}%
}}}}}
\put(1801,-2086){\makebox(0,0)[lb]{\smash{{\SetFigFont{10}{12.0}{\rmdefault}{\mddefault}{\updefault}{\color[rgb]{0,0,0}$\F$}%
}}}}
\put(1801,-1636){\makebox(0,0)[lb]{\smash{{\SetFigFont{10}{12.0}{\rmdefault}{\mddefault}{\updefault}{\color[rgb]{0,0,0}$\F$}%
}}}}
\put(1801,-1186){\makebox(0,0)[lb]{\smash{{\SetFigFont{10}{12.0}{\rmdefault}{\mddefault}{\updefault}{\color[rgb]{0,0,0}$\F$}%
}}}}
\put(1801,-736){\makebox(0,0)[lb]{\smash{{\SetFigFont{10}{12.0}{\rmdefault}{\mddefault}{\updefault}{\color[rgb]{0,0,0}$\F$}%
}}}}
\put(1801,-286){\makebox(0,0)[lb]{\smash{{\SetFigFont{10}{12.0}{\rmdefault}{\mddefault}{\updefault}{\color[rgb]{0,0,0}$\F$}%
}}}}
\put(1801,-2536){\makebox(0,0)[lb]{\smash{{\SetFigFont{10}{12.0}{\rmdefault}{\mddefault}{\updefault}{\color[rgb]{0,0,0}$0$}%
}}}}
\put(1801,-2986){\makebox(0,0)[lb]{\smash{{\SetFigFont{10}{12.0}{\rmdefault}{\mddefault}{\updefault}{\color[rgb]{0,0,0}$\F$}%
}}}}
\put(1276,-286){\makebox(0,0)[lb]{\smash{{\SetFigFont{10}{12.0}{\rmdefault}{\mddefault}{\updefault}{\color[rgb]{0,0,0}$U$}%
}}}}
\put(1276,-736){\makebox(0,0)[lb]{\smash{{\SetFigFont{10}{12.0}{\rmdefault}{\mddefault}{\updefault}{\color[rgb]{0,0,0}$U$}%
}}}}
\put(1276,-1186){\makebox(0,0)[lb]{\smash{{\SetFigFont{10}{12.0}{\rmdefault}{\mddefault}{\updefault}{\color[rgb]{0,0,0}$U$}%
}}}}
\put(1276,-1636){\makebox(0,0)[lb]{\smash{{\SetFigFont{10}{12.0}{\rmdefault}{\mddefault}{\updefault}{\color[rgb]{0,0,0}$U$}%
}}}}
\put(1276,-2536){\makebox(0,0)[lb]{\smash{{\SetFigFont{10}{12.0}{\rmdefault}{\mddefault}{\updefault}{\color[rgb]{0,0,0}$U$}%
}}}}
\put(3226,-1036){\makebox(0,0)[lb]{\smash{{\SetFigFont{10}{12.0}{\rmdefault}{\mddefault}{\updefault}{\color[rgb]{0,0,0}$HFI^+_{red}$}%
}}}}
\put(2176,-1336){\makebox(0,0)[lb]{\smash{{\SetFigFont{10}{12.0}{\rmdefault}{\mddefault}{\updefault}{\color[rgb]{0,0,0}$Q$}%
}}}}
\put(2176,-511){\makebox(0,0)[lb]{\smash{{\SetFigFont{10}{12.0}{\rmdefault}{\mddefault}{\updefault}{\color[rgb]{0,0,0}$Q$}%
}}}}
\put(2506,-2976){\makebox(0,0)[lb]{\smash{{\SetFigFont{10}{12.0}{\rmdefault}{\mddefault}{\updefault}{\color[rgb]{0,0,0}$\dl+1$}%
}}}}
\put(2521,-1626){\makebox(0,0)[lb]{\smash{{\SetFigFont{10}{12.0}{\rmdefault}{\mddefault}{\updefault}{\color[rgb]{0,0,0}$\du$}%
}}}}
\end{picture}%

 \end{center}
 
Here, $HFI^+_{red}$ is a finite dimensional space, and the infinite tower of $\mathbb{F}$'s is comprised of two sub-towers, connected by the action of $Q$. Each sub-tower is isomorphic to $\mathbb{F}[U,U^{-1}]/U$, so we call them $U$-towers. 
 
\begin{definition}
Let $\dl + 1$ and  $\du$ be the lowest degrees in the two $U$-towers of $HFI^+(Y, \mathfrak{s})$. These are the {\em involutive correction terms} of the pair $(Y, \mathfrak{s})$. 
\end{definition}

The involutive correction terms satisfy the properties
$$d\equiv \dl \equiv \du \ \  (\text{mod }2\Z), \ \ \  \underline{d} \leq d \leq \overline{d},$$
and (when applied to integral homology spheres) descend to maps
$$ \dl, \du: \Theta^3_{\Z} \to 2\Z \ \ \text{and} \ \  \dl, \du: \Theta^3_{\Z_2} \to \Q.$$

\bibliographystyle{amsalpha}

\end{document}